\documentclass[12pt,fleqn,a4paper, twoside]{article}
\usepackage{a4}
\usepackage{times}

\usepackage{amsmath}
\usepackage{array}
\usepackage{amssymb}
\usepackage{graphicx}
\usepackage[left=2.75cm,right=2.25cm,top=2.9cm,bottom=2.25cm]{geometry}

\def\cocoa{{\hbox{\rm C\kern-.13em o\kern-.07em C\kern-.13em o\kern-.15em A}}}

\newcommand{\nat}{\mbox{${\rm I\hspace{-.6 mm}N}$}}
\newtheorem{prop}{Proposition}[section] 
\newtheorem{set}[prop]{Setting} 
\newtheorem{rem}[prop]{Remark}
\newtheorem{thm}[prop]{Theorem}
\newtheorem{lemma}[prop]{Lemma}
\newtheorem{coro}[prop]{Corollary}
\newtheorem{defin}[prop]{Definition}

\newtheorem{remnot}[prop]{Remark-Notation}
 
\newtheorem{notat}[prop]{Notation}
\newtheorem{ex}[prop]{Example}

\newtheorem{(*)}[prop]{}
\newcommand{\m}{\mbox{${  m}$}}
\newcommand{\X}{\mbox{$ X\ $}}

\newcommand{\s}{\mbox{$\widetilde s$}}  
\newcommand{\n}{\mbox{$\#$}}
\newcommand{\af}{\mbox{\sf  {\it I}\hspace{-0.65 mm}A}}
\newcommand{\e}{\mbox{\sf  {$\ell$}}}

\newcommand{\pr}{\mbox{\sf {\it {\  I}\hspace{-0.6 mm}P}}}

\newcommand{\integ}{\mbox{${\sf Z\hspace{-1.3 mm}Z}$}}
\newcommand{\integsmall}{\mbox{${\rm\scriptstyle Z\hspace{-.9 mm}Z}$}}

\newcommand{\I}{\mbox{$\  \Longrightarrow \ $}}
\newcommand{\II}{\mbox{$\  \Longleftrightarrow \ $}}

\newcommand{\denu} {\end{enumerate}}

\newcommand{\0} {\raisebox{.5mm}{$ \scriptstyle{\bigcirc }$} }

\begin{document}






\begin{center}

\renewcommand{\thefootnote}{\fnsymbol{footnote}}
\addtocounter{footnote}{0}

{\LARGE On the minimum distance of AG codes, on \\[2mm]   Weierstrass semigroups and the smoothability of certain monomial curves  in 4-Space.}
\end{center}

\begin{center}

\renewcommand{\thefootnote}{\fnsymbol{footnote}}
\addtocounter{footnote}{0}
{\large
 Alessio Del Padrone,$^1$ Anna Oneto,$^2$   Grazia~Tamone$^3$ \footnote{A part of this work was done while the
last two authors were visiting the Department of Mathematics of the Indian Institute of Science, Bangalore, India.
They also    thank  Professor Dilip Patil for his   warm hospitality.} 
}  \\[2mm]
{\small $^{1,3}$Dima - University of Genova,
via Dodecaneso 35, I-16146 Genova,  Italy.}\\
{\small $^2$Diptem - University of Genova, Piazzale Kennedy, Pad.D 16129 Genova, Italia}  \\
{\small E-mail: {\tt $^1$delpadro@dima.unige.it,} {\tt $^2$oneto@diptem.unige.it,}  \ {\tt $^3$tamone@dima.unige.it}}
\end{center}

\vspace*{-5mm}

\noindent \rule[2pt]{\textwidth}{1pt}

\vspace*{-4mm}

\noindent \rule[0pt]{\textwidth}{1pt}

\vspace*{-5mm}

{\baselineskip8pt

\begin{abstract}

{\baselineskip5pt

\noindent In this paper we  treat several topics regarding numerical Weierstrass semigroups and the theory of Algebraic Geometric Codes associated to  a pair $(X,P)$, where $X$ is a projective curve defined over the algebraic closure of the finite field $\mathbb{ F}_q$ and $P $ is a  $\mathbb{F}_q$-rational point  of $X$. First  we show how to evaluate the  Feng-Rao Order Bound, which is a good estimation for the minimum distance of such codes. This bound is related to the classical Weierstrass semigroup of the curve $X$ at $P$.  Further we focus our attention on the question to recognize the Weierstrass semigroups over fields of characteristic 0. After  surveying  the main tools (deformations and  smoothability of monomial curves) we prove that the semigroups of embedding dimension four generated by an arithmetic sequence  are Weierstrass.  }
\end{abstract}
}

\vspace*{-5mm}

\noindent \rule[0pt]{\textwidth}{1pt}

{\baselineskip6pt
\vskip4pt

\noindent {\footnotesize\it Keywords~}{\scriptsize\sf :} {\scriptsize\sf
AG code, Order Bound, Numerical semigroup, Monomial curve, Deformation, Weierstrass semigroup. }
\smallskip

\noindent {\footnotesize\it {\rm 1991} Mathematics Subject Classification~}: {\scriptsize\sf Primary~:
14H55~; Secondary~: 14H37, 11G20, 94B27.}

}

\vspace*{-5mm}

\renewcommand{\thefootnote}{\arabic{footnote}}
\addtocounter{footnote}{0}

\setcounter{section}{-1}

\section{Introduction.}
The paper  is divided into two parts. In the first one we describe some bounds of the minimum distance of $AG$ codes, while in the second one we deal with the problem to  characterize the Weierstrass semigroups. 

In the first part $\mathbb{F}$ will denote the algebraic closure of  the finite field with $q$ elements $\mathbb{F}_q$; $X$ will be a smooth projective  algebraic curve  of genus $g$ defined over $\mathbb{ F}_q$.\\
 To a pair  $(X,P)$, where  $P\in X$ is a $\mathbb{F}_q$-rational point  can be associated a family of {\it Algebraic Geometric Codes} $C_i, \ i\in\nat$ and a 
 {\it numerical semigroup S}.   
 For $i$ large enough, the minimum distance $d(C_i)$ of such codes can be bounded by the  Feng-Rao order bound $d_{ord}(C_i)$  which depends only on the semigroup $S$ (see \cite{fr}). When $S$ is non-ordinary, it is called the Weierstrass semigroup of $X$ at $P$.
 Evaluations or estimates of the order bound are given by several authors,  either in general or in  particular cases (see, e.g., \cite{ba}, \cite{ot3}). In the first part of this paper we give a survey of these results and we state a conjecture (\ref{cong}) on the behaviour of the sequence $\{d_{ord}(C_i)\}_{i\in \nat}$, for $i>c+d-e-g$, where $c,d,e$ are suitable  integers associated to the semigroup $S$ (in    \cite{ot3}  this conjecture is proved in many cases).  

  According to the recalled relation with code theory, the classical study of Weierstrass semigroups is becoming relevant. In particular an interesting and  still open hard  question is  how to recognize   Weierstrass semigroups, i.e. those  semigroups associated to a smooth projective curve at a point $P$. This problem is approached in the second part under the simplifying assumption that $X$
is a smooth projective  algebraic curve  of genus $g$ defined over
an algebraically closed field $\mathbb{F}$ 
with  char$(\mathbb{F}) =0$.
    It is known that there exist non-Weierstrass semigroups: 
    the first example is due to Buchweitz, see \cite{bu}.
     A fundamental result on this topic has been proved by Pinkham in his Phd thesis \cite{pi}:  \\ \centerline{{\it \lq\lq$S$ is Weierstrass if and only if the monomial curve $X=\textrm{Spec} (\mathbb{F}[S])$ is smoothable".}}
    In some case it is known that a monomial curve is smoothable: 
see \cite{sc1} for  the   complete intersection case,
see \cite{sm} for  $X\subseteq \mathbb{A}^3$,
see \cite{k1} for $X\subseteq\mathbb{A}^4$ and Gorenstein, 
see \cite{k3} for semigroups of genus $g\leq 8$, 
see \cite{k2}, \cite{k4} for certain semigroups of embedding dimension $5$ or with $g=9$.

    In this paper we collect the main definitions and results on this question, further we illustrate the explicit algorithm to obtain a  deformation of a monomial $X$ with its $\mathbb{G}_m$ action and show several examples in a detailed way.
    Finally  we show that monomial curves in $\mathbb{A}^4$, generated by an arithmetic sequence are smoothable. It follows that every  semigroup $S$ of 
embedding dimension $4$ generated by an arithmetic sequence is Weierstrass. 
    
           \section{Weierstrass points and Weierstrass semigroups}
           
Let $\mathbb{F}$ denote an algebraically closed field.
Let $X$ be a smooth projective algebraic curve of genus $g$ defined over 
$\mathbb{F}$ with function field $\mathbb{F}(X)$, and     
let $P\in X$.
For each $k\in \nat$, let 
$$
{\cal L}(kP) \! =\!\{f\in\mathbb{F}(X)\setminus{0}\ |\ \textrm{div}(f)+kP\geq 0\}\cup\{0\}.
$$
This is clearly a vector subspace of 
$\mathbb{F}(X)$; we denote by $\lambda(kP)$ its dimension over $\mathbb{F}$. The following are well-known facts:

$\lambda(kP)=\textrm{dim}_{\mathbb{F}}({\cal L}(kP))\in\nat$,\quad
$\lambda((k-1)P)\leq \lambda(kP)\leq \lambda((k-1)P)+1$ for each $k>1$,
and by Riemann-Roch Theorem $\lambda(nP)=n-g+1$ for each $n\geq 2g-1$.

Hence the set $H(P):=\{k\in\nat_+|\ \lambda((k-1)P)=\lambda(kP)\}$, 
of {\it gaps at $P$}, is a proper subset of $\{1,2,\dots,2g-1\}$ and it has exactly $g$ elements. Moreover it is easy to see that its complement $S(P):=\nat\setminus H(P)$, the {\it set of non-gaps at $P$}, is a numerical semigroup.\\

Recall that a semigroup $S$ is called {\it ordinary} if it is of the form 
$S=\{0,e,e+1\rightarrow\}$ for some $e>0$ (note that its {\it genus}, also called $\delta$, is exactly $e-1$). 

\begin{defin}
A {\it Weierstrass point} of $X$ is a point $P$ such that $H(P)\neq \{1,\dots, g\}$. A semigroup $S$ is called {\it Weierstrass} (over $\mathbb{F}$) if
there exists a smooth projective algebraic curve $X$ (defined over $\mathbb{F}$) and a Weierstrass point $P$ such that $S=S(P)$. 
\end{defin}

See for more details, e.g., \cite[Exercise A.4.14]{HiSi00} or \cite{dc}.

\begin{rem}\label{buco}
{\rm 
Let $P\in X$, then by Riemann-Roch Theorem
\begin{enumerate}
\item
$n\in S(P)$ $\II$ there exists $f\in \mathbb{F}(X)$ such that $(f)_\infty=nP$, i.e. $\textrm{ord}_P(f)=-n$.
\item 
$n\in H(P)$ $\II$ there exists a regular differential form $\omega$ with $\textrm{ord}_P(\omega)=n-1$ (because by  Riemann-Roch theorem: $  \lambda(K-(n-1)P)>0,$    for each $  \ gap \ n\in H  $, where $K$ denotes a canonical divisor).

\item $P$ is a Weierstrass point $\II$ $\lambda(gP)\geq 2$ $\II$ there exists a regular differential form $\omega$ with $\textrm{ord}_P(\omega)\geq g$. In particular, it follows immediately that
if $X$ has a Weierstrass point then $g\geq 2$.

\item By the previous point, the presence of a Weierstrass point on an algebraic curve  of genus $g$ ensures the existence of a morphism of degree not exceeding $g$ from the curve  onto the projective line: pick the morphism associated to the linear system $|iP|$ with any $i$ such that $\lambda(iP)=2$ and $i\leq g$.

\end{enumerate}
}
\end{rem}

\subsection {On the number of Weierstrass points on a curve.}
For a smooth curve $X$ let $W$ denote the set of Weierstrass points of $X$. We know that
\begin{enumerate} 
\item If $g\leq 1$ the set  $ W $ is empty.
\item Case $X$   {\it hyperelliptic}. 
A hyperelliptic curve is an algebraic curve which admits a double cover over $\pr^1$. These curves are among the simplest algebraic curves: they are all birationally equivalent to curves given by an equation of the form $y^2=f(x)$ in the affine plane, where $f(x)$ is a polynomial of degree  $>4$ with  distinct roots,
 and the degree of $f(x)$ is either twice the genus of the curve plus 2, or twice the genus of the curve plus one.

If a double cover  exists, then it is the unique double cover  and it is called the
\lq\lq hyperelliptic double cover".
In algebraic geometry the Riemann-Hurwitz formula,  states that if $X$, $ X'$ are smooth algebraic curves, and $\Phi: X \longrightarrow  X'$ is a finite map of degree $d$ then the number of branch points of $\Phi$, denoted by $N$, is given by
$$2 g(X)-2 =2d\ \!(g(X')-1)+N$$
By the Riemann-Hurwitz formula the hyperelliptic double cover has $X'=\pr^1$, hence has exactly $2 g+2$ branch points. For each branch point $P$  we have $\lambda(2P)=2$, hence these points are all Weierstrass points; for each of them  there exists a function $f$ with a double pole at P only. Its powers have poles of order $4, 6$, and so on. Therefore at $P$   the gap sequence is
\quad $1, 3, 5, ..., 2g -1$ \  and $\lambda(kP)=2k$, we conclude that the Weierstrass points of $X$ are exactly the $2g+2$ branch points of the hyperelliptic double cover.
\item For algebraic curves of genus $g$  there always exist at least  \ \ $2g+2\ \ $ Weierstrass points and only the hyperelliptic curves of genus $g$ have exactly  $2g+2$ Weierstrass points. 
\item The upper bound on the number of Weierstrass points is $  \ g^3-g.\ \ $ 
\item  \cite{mac}  For each $g\geq 3$ there exist compact Riemann surfaces of genus $g$ with at least two Weierstrass point with different gap sequences. 
\end{enumerate}

           \section{Algebraic-geometric codes}
   
Let now $\mathbb{F}$ denote an algebraic closure of the finite field with $p$ elements 
$\mathbb{F}_p$, $p$ prime. Let $X$ be a smooth projective algebraic curve of genus $g$ defined over $\mathbb{F}_q$, $q=p^r$ for some $r\in\nat_+$, with function field 
$\mathbb{F}(X)$.           
Let $P\in X$ be an $\mathbb{F}_q$-rational point:
a family of {\it codes} and a {\it numerical semigroup} can be associated to $(X,P)$ as follows. 
For each $k\in \nat$, we consider the vector subspace of $\mathbb{F}_q(X)$
defined as
$${\cal L}_{\mathbb{F}_q}(kP) \! =\!\{f\in\mathbb{F}_q(X)\setminus{0}\ |\ \textrm{div}(f)+kP\geq 0\}\cup\{0\},
$$
it can be shown that \;
$\lambda(kP)=\textrm{dim}_{\mathbb{F}}({\cal L}(kP))=\textrm{dim}_{\mathbb{F}_q}({\cal L}_{\mathbb{F}_q}(kP))$ (\cite[Proposition A.2.2.10.]{HiSi00}).\\

We now recall the definitions of the AG codes associated to the pair $(X,P)$.
Choose $P_1,...,P_n$ distinct $\mathbb{F}_q$-rational points on $X$ such that $P_j\neq P $ for each $j$, and consider the $\mathbb{F}_q$-linear map
$$
\Phi_k\colon  {\cal L}_{\mathbb{F}_q}(kP) \longrightarrow    {\mathbb{F}_{q}}^{ n} \quad {\rm as }\quad  \Phi_{k}(f) = (f(P_1),..., f(P_n)).
$$
\begin{defin}
The family of {\it one-point AG   codes  of order  $n$ } is defined as
$$
C_k:=(Im\ \Phi_k)^{\perp}=\{x\in {\mathbb{F}_{q}}^{ n}\mid <x, \Phi_k(f)>=0 
\textrm{ for all }f\in {\cal L}_{\mathbb{F}_q}(kP)\},
$$
where $<x,y>:=x_1y_1+\dots+x_ny_n$ for each $x,y\in {\mathbb{F}_{q}}^{ n}$.
\end{defin}

A good estimate of  the {\it minimum distance} $d(C_k)$  of an AG code   is the  Feng-Rao order bound      ${ d_{ORD}(C_k)}$ which depends only on the semigroup $S=S(P)$. Let us fix the following notation
$$
S=\{s_0=0,s_1,..., s_j,...\}\neq \nat 
$$ 
with $s_i<s_j$ if $i<j$.
\begin{defin}
   For $s_j\in S$, let  $\left[\begin{array}{lcllll}N(s_j)  &:=&   \{(s_h,s_k)\in S^2\ |\  s_j=s_h+s_k\} \vspace{0.2cm}\\
     \nu(s_j)\!&:=&   \n N(s_j)\end{array}\right.$.    \vspace{ 0.2cm}
  \\ The    Feng-Rao  order   bound  of the code $C_k $ is      $ \begin{array}{lcllll}  { d_{ORD}(C_k)}\!& :=&\textrm{min} \{\nu(s_j)\ |\ j>k\}\leq d(C_k).\end{array}  $
  \end{defin}
If $S $ is ordinary, that is $S=\{s_{0}=0,s_{1}=g+1,s_{2}=g+2\rightarrow\}$, the sequence  $\{\nu(s_j),\ j\in\nat\}$ is non-decreasing and so  
$$
d_{ORD}(C_k)= \nu(s_{k+1}) \ \ for\     \ k\geq 0.
$$
In the other cases, it is known that   there exists  ${\it m} \ \!   \in \nat_+ $ such that
 $$ \nu(s_m) > \nu(s_{m+1}) \ \  and  \ \ \nu(s_{m+i})\leq   \nu(s_{m+i+1})\ \ \forall \ i \geq 1 .  $$ 
Then: $ d_{ORD}(C_k)= \nu(s_{k+1}) $ for each code $C_k$ with $k\geq m$.\\
  
\subsection{ Methods for the evaluation of $s_m$.}

Our goal is to find $ s_m $  for a given semigroup S; to this end it is useful
to consider the elements of $S$  $\lq\lq$near"  the conductor. 

\begin{notat}\label{not0} 
We shall refer to a numerical semigroup $S$,   with finite complement in $\nat$\ 
 $$S=\{0=s_0,s_1,..., s_j,...\}\neq \nat $$ 
where $s_i< s_k $, \ if \ \ $i < k $. 
Further we denote: 
$$
\begin{array} {lll}
\textrm{embdim}(S)\ =\  minimal \ number\  of\  generators\  of\  S
\vspace{0.1cm} \\ 
e=s_1= \textrm{min}\{s\in S\mid s\neq 0\},\ the \ multiplicity 
\vspace{0.1cm} \\ 
 c= \textrm{min}  \{   r\in S\ |\  r+\nat\subseteq S \},    \  the \     conductor\vspace{0.1cm} \\  
 d=    \textrm{max}\{s_i\in S \ |\ s_i<c\},   \   the \ dominant\ \vspace{0.1cm}\\ 

c^\prime= \textrm{max}\{s_i\in S \ |\ s_i\leq d\ \  and\ \ s_i-1\notin S\},\       the \   subconductor  \vspace{0.1cm} \\ 

  d^\prime=    the\ greatest\ element \ in \ S\ preceding\  c  ',  \    when  \ \ c'>0\vspace{0.1cm} \\
 
 \ell= c-1-d,  \   the\ number\ of\   gaps\ of\ S\ greater\ than  \ d\vspace{0.1cm}\\  
 \s= \textrm{max}\{s\in S,\ |\ s\leq d,\ \   s-\e\notin S\}.
\end{array} 
$$
This means that $S$ has the following shape  (thinking of it as embedded in $\nat$, where $*$ means a \lq\lq gap" of $S$) 
$$
\begin{array}{ccccccc}
 &\scriptstyle{e-1}\ {\rm \scriptstyle{gaps}}&&
 \scriptstyle{c^\prime-d^\prime-1}\ {\rm \scriptstyle{gaps}}&&\ell \ {\rm \scriptstyle{  gaps}}&\vspace{-0.1cm}\\
  S=\{ 0, &*\dots *& \ e,   \dots   ,\  d^\prime, & *\dots* & c^\prime\ \longleftrightarrow\ d, &  *\dots *& c\rightarrow  \} 
  \end{array}$$
\end{notat}
 A  semigroup $S$ is called  { \it acute\ }   if  either \  $S$ \ is ordinary, \ or    \  \  $c,d,c  ',d  '$ \ \ satisfy    \    \ $c-d\leq
c  '-d  ' $ (see \cite{ba}). 
If $S$ is non-ordinary, it can be seen that: \\
\centerline{$S$ acute $\I$  $c'\leq \s\leq d$.} 
 
 \begin{ex} $S=\{0,8_e, 12_{d'},14_{c'},15,16_d, 20_c\rightarrow\} $ \ has $(\ell=3,\ \s=14,\ c'-d'=2<c-d=4$, $S$ non-acute$ )$.
\end{ex}  

  From now on, $S$ will be non-ordinary.
In order to evaluate $s_m$ we    study the difference   $\nu(s_{i+1})-\nu(s_i)$   for  $s_i\in S $. To this end, it 
is $\lq\lq$natural"  to consider the following partition of  {$N(s_i) =\{(s_j,s_k)\in S^2\ |\  s_i=s_j+s_k\}$:}
  $$N(s_i) = A(s_i)\cup  B(s_i) \cup C(s_i)\cup  D(s_i) $$  $$\begin{array}{lcllll}   A\ \! (s_i)&:=&\{(x,y), (y,x)   \in N(s_i) \ |\    x<c  ',\     c  '\leq   y\leq d \} 
   \vspace{0.2cm}\\

B\ \!(s_i)&:=&\{(x,y)\in N(s_i)\ |(x,y)\in[c  ',d]^2\ \}   \vspace{0.2cm}\\

   C\ \!(s_i) \!&:=&   \{(x,y)\in  N(s_i)\ |\   \   x\leq d',\ y\leq d'\} 
 \vspace{0.2cm}\\
           D\ \!(s_i)\!  &:=&    \{ (x,y),\ (y,x)\in N(s_i)\ |\   \  x\geq c,\ x\geq y \}. 
  \end{array} $$
    \begin{ex} $S=\{0,8_e, 12_{d'}, 14_{c'},15,16_d, 20_c\rightarrow\} $.   For $i=16, \ { s_i=30:}$  \\    $A(s_i)=C(s_i)=\emptyset$,\ \ $B(s_i)=\{(14,16),(15,15),(16,14)\}, $  \ \\ $D(s_i)=\{(0,30),(8,22),(30,0),(22,8)\}.$ \\
  For $ { s_6=20:}$  \ \ $A(s_6)=B(s_6)=\emptyset,$  $\ C(s_6)=\{(8,12),(12,8)\}\ ,\ \ D(s_6)=\{(0,20),(20,0)\}.$

\end{ex}

  \begin{set} \vspace{-0.2cm} $$\begin{array}{lcllll}      
  \alpha(s_i)& :=&\n A(s_{i+1})-\n A(s_i)  \vspace{0.2cm}\\

  \beta(s_i)& :=&\n B(s_{i+1})-\n B(s_i)  \vspace{0.2cm}\\

  \gamma(s_i)\!& :=&   \n C(s_{i+1})- \n C(s_i) \vspace{0.2cm}\\
               \delta(s_i)&: =&\n D(s_{i+1})-\n D(s_i). 
  \end{array}$$ 
 \end{set}
Therefore: \ \ { $\nu(s_{i+1})-\nu(s_i)=\alpha(s_i)+\beta(s_i)+ \gamma(s_i)+\delta(s_i)$.}
\begin{lemma}\label{varie} $($see {\rm \cite{ot3})}  \begin{enumerate}
\item $\alpha(s_i)\in \{-2,0,2\}$ \ \ and \ \ $\alpha(s_i)=0$, \ if \ $s_i>d'+d$.
\item $\beta(s_i)\in \{-1,0,1\}$ \ \ and \ \ $\beta(s_i)=0$, \ if \ $s_i>2d$.
\item $\gamma(s_i)$ {\rm is difficult to evaluate if $s_i< 2d'$, trivial otherwise:
 \item[] in fact \ $\gamma(2d')=-1$ \ \  and \ \ $\gamma(s_i)=0$, \ if
 $\ s_i>2d'$.}
 \item If   $s_i\geq 2c$, then $\delta(s_i)= 1.$ \
\item[]
If $s_i<2c $ and $s_{i+1}\in S$, then $\delta(s_i)\in \{0,2\}$ \ and 
 $\ \ \delta(s_i)=0 \II s_{i+1}-c \notin S$. 
\item $ s_m\leq 2d$.  
{\rm(In fact by (1)-(4),  if $   \ s_i\geq 2d+1$ then $\alpha=\beta=\gamma=0$ and so
$\nu(s_{i+1})-\nu(s_i)= \delta(s_i)\geq 0.$}
\denu
 \end{lemma}
 By (\ref{varie}.5),  from now one has to consider only  elements $\ s_i\leq 2d$, in order to   find the greatest $s_i\in S$ such that $\nu(s_{i+1})<\nu(s_i)$. Assume $\ s_i+1\in S$.  
\begin{rem}  \begin{enumerate}
\item If $s_i=\s+d$, then 
  $s_{i+1}-c= \s-\e\notin S$ \ (by definition),  and so  \  $s_i=\s+d$     is the greatest element   satisfying   $\delta(s_i)=0$.\\ 
  {\rm (For this reason \  $\s+d$ \ is a \lq\lq good candidate" for $ \ s_m$)}.
\item If $s_i\geq 2d'$ we know that $\gamma(s_i)\leq 0$ and easily one can see when $\nu(s_{i+1})<\nu(s_i).$ 

\item If $s_i< 2d'$ we can write $\nu(s_{i+1})-\nu(s_i)$ in function of $\gamma(s_i)$: it depends also on the facts:  \begin{enumerate}
\item[]   $s_{i+1}-c$   $\in S$ \  or \  $\notin S$,
\item[]   $s_i-d$  \ \ \ $\in S$ \  or \  $\notin S$,
\item[]
  $s_{i+1}-c'$  $\in S$ \  or \  $\notin S$.\denu
  \denu 
    
  \end{rem}
  In \cite{ot3} the results  on the position of $s_m$   are explained by means of several tables. For example we show for $s_i<2d'$ how the difference $\eta(s_i):=\nu(s_{i+1})-\nu(s_i)$ depends on the value of $\gamma:=\gamma(s_i)$. \quad In the following table \quad $\times$ means $\in S$ \ and $\ \0$ \ means \ $\notin S$.
Assume \  $s_{i}\leq 2d  '-1  $ Then: \vspace{-0.2cm} 
$$
\begin{array}{c|c|c|r|r|r||c}	 	
s_{i+1}-\!c	 &	s_{i}-d	& s_{i+1}-\!c  ' &
\alpha	&	\ \beta	&	\ 	\delta& \quad\eta(s_i)\quad\\
\hline
\hline
\0	&	 	\times	&	\0	&	-2	&	0	&	 	0	&	\gamma-2  	\\
\0	&	 	\times	&	\times	&	0	&	0	&	 	0	&	\gamma	\\
\0	& 	\0	&	\0	&	0	&	0	& 	0&\gamma  	\\
	\times	&	 	\times	&	\0	&	-2	&	0	&	 	2	&	\gamma	\\
\0	&	 	\0	&	\times	&	2	&	0	 	&	0	&	\gamma+2	\\

	\times&	 	\0	&	\0	&	0	&	0	&	 	2	&\gamma+	2	\\

	\times&	 	\times	&	\times	&	0	&	0	&	 	2	&	\gamma+2	\\
	\times&	 	\0	&	\times	&	2	&	0	&	 	2	&\gamma+	4	\\
\end{array}
$$

 Recall: $\gamma(s_i)$ concerns pairs $(x,y)\in\nat(s_i)\cap[0,d']^2$.

\subsection{ Evaluation or bounds for $s_m$.}
 \begin{thm} \label{ev}$(See\ \cite{ot3})$ \ With  setting $(\ref{not0})$ we have:
 \begin{enumerate} \item   If \ \ $\s<2d'-d$,\ \ then \ \ $s_m\leq 2d'.$  
  \item[]  If, moreover, $\ [\s+2,d']\cap\nat\subseteq S$,\ \ then\ \ $s_m=\s+d$.
  \item  If \ \ $\s\geq 2d'-d$,\ \ then \ \ $s_m\leq \s+d$.\ \   
  \item[]   More precisely:
    \begin{enumerate}\item If \ \ $\s\geq  d'+c'-d$,\ \ then \ \ $s_m=\s+d.$  \item If \ \ $\s= 2d'-d$,\hspace{0.6cm}\ \ then \ \ $s_m=\s+d.$    

 \item If \ \ $  2d'-d< \s< d'+c'-d$,\ \ we  can give upper and lower bounds for $s_m$ under additional assumptions. \ \ In particular: 
    \item[]if \ \ $  [d'-\e,d']\cap\nat\subseteq S$,\ \  then \  \  $  c+d-e\leq\s+d'-\e+1\leq s_m\leq 2d'$.    \    \ \  \ \    
 \end{enumerate}    Case (a) is satisfied e.g. if $ \ d-2\leq \s\leq d; $ or if ${\ c'\leq \s\leq d}$, in particular if  $S$ is acute.
     \end{enumerate}
     \end{thm}

   \begin{ex}{\rm   \begin{enumerate}     
    
\item  $S=\{0,25_e,26,28, 30,31_{d'},33_d,39_c \rightarrow\}$
\\    $(\s=28, \ \s<2d'-d,\ [\s+2,d']\cap\nat\subseteq S, \ {s_m=\s+d })$. 
   \item   $S=\{0,7_{e=d'},13_{c'},14,15,16,17_d,20_c \rightarrow\} $\\
   $(S$ is acute, $\e= 2$, $\s=14$, \ $c'\leq \s\leq d$, $\ \s>  d'+c'-d)$.

 \item  $S=\{0,20_e,21,26,27_{d'},32_d,39_c \rightarrow\}$\\    $(\s=21<2d'-d,\ {s_m=2d'=54>\s+d} )$. 
  \item  $S=\{0,10_e,20,22,23_{d'}, 26_d,30_c \rightarrow\}$\\    $(2d'-d<\s=22<d'+c'-d ,   \ \ {s_m =46<\s+d})$. 
 \denu}
 \end{ex}

   \subsection { Conjecture and particular cases.}
   We believe the following fact is   true for each semigroup.\vspace{0.2cm}\\
   \centerline{ { \bf Conjecture:} \label{cong}\qquad $s_m\geq \ c+d-e\qquad(*)$\hspace{3cm}  }\vspace{0.2cm}
We   proved in $\cite{ot3}$ \ that $(*)$ holds in several cases, in particular 
 {\it  \begin{enumerate}
\item If either \ $(s_m\geq \s+d)$ \ or \ $(s_m\geq 2d'$ and $\ \s<d')$.
 \item If \ $2d'-d<\s<d'+c'-d$ \ and \ $[d'-\ell,d']\cap\nat\subseteq S$ $(\ref{ev}.2c)$.  
   \item When \ $\e=2$, \ or \ $\e=3$ (here  we   calculate $s_m$ exactly).  
   
    \item If \ $\tau \leq 7$ \\ $($where $\ \tau:= \n\{x\in \nat\setminus S\ |\ x+(S\setminus\{0\})\subseteq S\}$ is the {\it Cohen-Macaulay type of \ $S$}$)$.    
   
  \item If\  $e\leq 8  \ \ ($by $(4)$, since \ $\tau\leq e-1)$.
  
    \item If \ $S$ is generated by a  {\bf generalized arithmetic sequence}  $($i.e. $S=<m_0,m_1,...,m_p>$ where $m_i=am_0+id$, \ for some $a\geq 1, d\geq 1$)$,\ $ then $s_m=\s+d$ and so $(*)$ holds.
  
    \item  If \ $S$ is generated by an  {\bf almost arithmetic sequence}      $($i.e.   $S=<m_0,m_1,...,m_k,\ n>$, where $\ m_0,m_1,...,m_k\ $ is an arithmetic sequence$)$  and  $\ { \textrm{embdim}(S)\leq 5}$, \\   then \ \ $s_m\geq c+d-e$.    
	     \denu}
  
 \section{Weierstrass Semigroups.} 
  
In this section we deal with the following \\  
\centerline{ {\bf Question :}   Which numerical semigroups are Weierstrass?} \vspace{0.2cm}
    \noindent The problem to find conditions in order that a semigroup is Weierstrass seems to be very hard: there are  only   partial answers in several directions. Most of them are in  $characteristic \ 0$, so we fix the following  

\begin{set} From now on we assume that $\mathbb{ F}$ is algebraically closed with $char(\mathbb{F})=0$.
 \end{set}
 We know that there exist non-Weierstrass semigroups:
 the first  example is due to an idea of Buchweitz: 
 \begin{ex}{\rm (See \  \cite{bu}) \   Let $S=<13,14,15,16,17,18,20,22,23>$, \ with  \ $g=16$,  $c=26, \\  H=\!\nat\!\setminus\! S=\{1,...,12,19,21,24,25\}$.  
 
$S$ cannot be Weierstrass.  In fact assume that  there exist a curve $X$ and a point $P\in X$  such that $S=S(P)$. Then, by Remark \ref{buco}, $X$ would have regular differentials $\omega_i$ vanishing at $P$ to orders  
$i=\textrm{ord}_{P}(\omega_i) $ with $i\in \{0,1,2,\dots,10,11,18,20,23,24\}$. \\ 
 Hence, taking suitable (tensor) products of the differential forms above,
$X$ would have also at least $46$ linearly independent \lq\lq quadratic" differentials 
vanishing to every order \\ $\in \{0,...,35,36,38,40,41,42,43,44,46,47,48\}$ at $P$. 
This implies that $\lambda(2K)\!\geq\!46$, a contradiction since, by Riemann-Roch it is
$\lambda(2K)=3g-3=45.$ }
  \end{ex}

   There are generalizations of this idea due to   Kim \cite{ki} and Komeda \cite{k1}:   
  
 \begin{prop}\label{kom}{\rm \cite{k1}}
   For a semigroup  $S$    of genus $g$,   let 
   $ \ \nat\!\setminus \! S= \{h_1,\dots,h_g\}$   and let \\
 \centerline  {  $H_m:= \{h_{i_1}+\cdots+ h_{i_m}\ |\ 1\leq  i_j\leq g\}$,\ \ $m\geq 2$.}
    \\ If \ $S$  \ is   Weierstrass, then 
 {$ \n \ H_m\leq (2m-1)(g-1)\ \ for\ each \ m\geq 2 \qquad(**) $}   \end{prop}
 { {\it {\it Proof.}} } If $S$ is the Weierstrass semigroup of a curve $X$ at $P$, then 
 $X$ has regular differentials vanishing to order $h_i-1$,  $\forall   i=1,..,g$. In fact let  $K$ be a canonical divisor  (in particular $\textrm{deg}(K)=2g-2$): for each $h_i\in \nat\setminus S,$ \ $\lambda(h_iP)=\lambda((h_i-1)P)$ therefore\vspace{0.1cm}
 by  Riemann-Roch $$\ \lambda(K -(h_i-1)P)>0.$$    
 It follows $\ \lambda(mK)\geq \n H_m,\ \forall m\geq 2$,  because    $\forall y_j\in H_m$,   ${\cal L}(mK) $ contains a  $m$-differential   vanishing to order $(y_j\!-\!m)$ at $P$. 
 Now it suffices to recall   that, again by  Riemann-Roch,  $\lambda (mK)=(2m-1)(g-1)$.\quad 
 $\diamond$  
 
 \begin{rem} The conditions $(**)$ of \ (\ref{kom}) \ are    satisfied for each $m\geq 2$, if \ $\ 2c< 3g$.
 \end{rem}
 {\it Proof.}
 Since  $\nat\setminus S\subseteq [1,c-1]\cap\nat$ we get $\n(H_m)\leq m(c-1)$, then the  inequality $(**)$ in \ (\ref{kom}) is surely  satisfied if $m(c-1)\leq(2m-1)(g-1)$ for each $m\geq 2$. \ 
 This condition is equivalent to  \quad $mc \leq (2m-1)g-(m-1)$: for 
   $m=2$, get   \  $2c\leq 3g-1$, i.e. $c\leq 3g/2-1/2$. Now assume $2c\leq 3g-1$,  and so    $mc\leq 3mg/2-m/2\  \forall m>0: $   one can easily see that the inequality $3mg/2-m/2\leq (2m-1)g-(m-1)\ $ holds   $\forall m\geq 2,\  \forall g>0$.\quad $\diamond$  
  
\begin{rem} \ In Buchweitz's  example,  $m=2$,  $g=16$, $\n H_2=46>3g-3$.   Note that for  $m=2$, the genus $g=16$ is the $\lq\lq$minimum"  example:  in fact 
Komeda and Tsuyumine found by a direct computation that \\
\centerline{{ for each  $2\leq g\leq 15\ $ we have $\ \n H_2\leq 3g-3$.}}
\end{rem}
F.Torres found a method to construct  symmetric non-Weierstrass semigroups (of large genus):
 \begin{ex} $(See\ {\rm\cite{to}})$
  {\rm Let $S'$ be a non-Weierstrass semigroup of genus $\gamma$, and let $g\in \nat$, \  $g\geq 6\ \!\gamma+4. $ \ Then the following symmetric   semigroup: 
 \\ 
 \centerline{  $S=\{2s\ |\ s\in S'\}\cup\{2g-1-2t\ |\ t\in\integ\setminus S'\}$}  
is  }   non-Weierstrass.
\end{ex}
This fact is true since we have:  
 \begin{prop} {\rm\cite[Scholium 3.5]{to}} Assume that a semigroup $S$ of genus $g\geq 6\ \!\gamma+4, $ is 
    $\gamma$-hyperelliptic,  i.e. satisfies
\begin{enumerate}
 \item  the first $\gamma$ elements $m_1,...,m_{\gamma}\in S,\ m_i>0$, \    are even;
\item $m_{\gamma}=4\gamma$;
\item  $4\gamma+2\in S$.
\denu Then:  
$S$ Weierstrass $\I S'  :=\left\{0,\displaystyle{\frac{m_1}{2},..., \frac{m_{\gamma}}{2}}\right\}\cup_{i\in\nat}\{2\gamma+i\} $ is Weierstrass.
 \end{prop}
  
\begin{ex}
{\rm\cite{o}  The possible Weierstrass semigroups for a plane smooth projective quintic (hence of genus $6$) are of the following types: 
  \begin{enumerate}
\item[]\hspace{1.2cm} $S_1=<4,5> $
\item[]\hspace{1.2cm} $S_2=<4,7,10,13> $
\item[]\hspace{1.2cm} $S_k=\{0,6\rightarrow\} \setminus \{k\}\ with\ 6\leq k\leq 11.$
\denu 
Note: $S_1,\ S_2$ and $S_k$ for $k=6,11$ are  semigroups   generated by an arithmetic sequence,  and $S_k$ is   generated by an almost arithmetic sequence for $k=10$.}
\end{ex}
\subsection{  Deformations   and $T^1({{\cal O}_X }_{\!,O}).$ }
The next theorem due to  Pinkham (thesis)  is fundamental to approach our question.
\begin{thm}\label{fond}{\rm \cite{pi}} \ 
 Let $S$ be a numerical semigroup and let \ $X= \textrm{Spec}(   \mathbb{ F}[S])$   be the monomial curve associated to $S$. Then:\\
\centerline{ { S is Weierstrass if and only if X is smoothable}}
\end{thm}
We want to recall the main tools of the theory.
Recall that the field $ \mathbb{ F}$  is algebrically closed with $\textrm{char}(\mathbb{F})=0$.
 We collect here the most important results and definitions on deformations of algebraic varieties.
\begin{defin} A deformation  $\pi:Y\longrightarrow \Sigma$ of a variety \X is a cartesian diagram
$$\begin{array}{cccc} \X&\hookrightarrow &Y\\
\Big\downarrow&&\ \ \Big\downarrow \pi\\
\{0\}&\hookrightarrow &\Sigma\end{array}$$
where $\pi$ is a flat morphism.\\
A deformation $\pi:Y\longrightarrow \Sigma$ of \X is said to be versal if any deformation $\pi':Y'\longrightarrow \Sigma'$ of \X is isomorphic to a deformation obtained from $\pi$ by a base change $h :\Sigma'\longrightarrow \Sigma $:  
$$\begin{array}{cccc}&_{_{\displaystyle{pr_1}}}\\
 {  Y'}= Y\times \Sigma'&\longrightarrow &Y\\ 
 \\
\ \ \ \ \pi'\Big\downarrow&&\ \ \Big\downarrow \pi\\
\\
\quad \ \ \Sigma'&\longrightarrow &\Sigma\\
&^{^{\displaystyle{h}}}
\end{array}$$
When $\Sigma= \textrm{Spec}\ k[\varepsilon]/(\varepsilon ^2)$  we say that the deformation is infinitesimal.\\
Finally we say that a deformation is trivial if $\ \ Y\simeq \X\times \Sigma$.
\end{defin}
\begin{defin} A variety $X$ is smoothable if there exists a deformation $Y$  of $X$ having smooth generic fibre.
\end{defin}
For a   survey on  deformations  we refer to {\rm \cite{st}}.
We recall the main theorem
  \begin{thm}\label{pis}{\rm\cite{pi}}\ 
  
{If $X$ is affine  variety and has an isolated singularity, then there exists a versal deformation $Y$ of $X$. \ 
Further, if $X$   has a   $ \mathbb{G}_m$-action, then there exists a    $\mathbb{G}_m$-action on $Y$ extending the action on $X$.}
\end{thm}
\begin{coro} Let  $X= \textrm{Spec}(   \mathbb{ F}[S])$, $S$ a numerical semigroup. Then  $X$ has a versal deformation $Y$ compatible with the well-known $   \mathbb{G}_m$-action.
\end{coro}  
  \begin{notat}\label{GH} Let $S=<n_0,...,n_k>$  be a semigroup,  $  P=   \mathbb{ F}[x_0,...,x_k]$,   $\textrm{weight}(x_i):=n_i$ \ \ $(0\leq i\leq k)$,  and  let $$B:=   \mathbb{ F}[S]=   \mathbb{ F}[t^{n_0},t^{n_1},...,t^{n_k}]  =P/I,$$  \ where \ \ $I=(f_1,...,f_q) ,\ \ f_i$     homogeneous binomials, \
    {  $ d_i:=\textrm{deg}(f_i),\ \forall i=1,..., q$.  Further let:
    \begin{enumerate}
    \item[]   $f= \begin{pmatrix}f_1\\ \vdots\\f_q\end{pmatrix}\in P^{q} $ 
  \item[] $  G_{\ell}:=\{i\in\{0,\dots,k\} \ | \ n_i+\ell\notin S\} $
    \item[] $  H_{\ell}:=\{d_k,\ k\in\{1,\dots, q\} \ | \ d_k+\ell\notin S\}  \ (\ell\in \integ).$
    \denu}
\end{notat}
In order to construct deformations for the curve $X$ we need  the $B$-$module$ $T^1_B$.
\  Let   $\Omega_{P/   \mathbb{ F}}$ be the $P$-$module$ of 1-$differential\ forms$,
   then \ 
    {$ Hom_B(\Omega_{P/   \mathbb{ F}}\otimes B,B)$ }
  is a free $B$-$module$ generated by the partial derivatives \  \ $<{\frac{\partial \  }{\partial x_{0}}},..., {\frac{\partial \  }{\partial x_{k}}}> $:
  \begin{defin} \label{T1}  
 Consider the map 
    $$\begin{array}{ccccc} \varphi: & Hom_B(\Omega_{P/   \mathbb{ F}}\otimes B,B) & \longrightarrow & Hom_B(I/I^2,\ B)  \\
 &\hspace{3cm}  {\frac{\partial \  }{\partial x_{i}}}&\mapsto & g\ :\   \ g({f)= (\frac{\partial f  }{\partial x_{i}}})&(\textrm{mod}\ I)
  \end{array}$$
   We define the $B$-module $T^1_B$ as   \ {  $T^1_B:= Coker \ \varphi$ }. 
    \end{defin}

 Let{  $f$ be as in (\ref{GH}):  we shall   identify a map $g\in  Hom_B(I/I^2,\ B)$ with the column vector  $ (h_j)_{j=1,... q} :=(g({ f}  ))  $ of its image $\textrm{mod}\ I$.
\begin{remnot} Let $J_0$ be the jacobian matrix of $\textrm{deg} \ 0$-$derivatives:$  $J_0=\left(\begin{array}{c} x_i\displaystyle{\frac{\partial f_j }{\partial x_{i}}}\end{array}\right)$. Then  $$J_0\equiv \left(\begin{array}{llll}
  t^{d_1} & 0&\ \dots \ &0 \\
 0& t^{d_2}&\ \dots  \ &0  \\
 &\dots& \dots\\
 &\dots& \dots\\
0 &0 &\ \dots\ &t^{d_q} 
\end{array}\right) J_0(1)\quad (\textrm{mod}\ I).$$
where   $J_0(1)$ \  is the evaluation of  $J_0$ at the regular point $Q(1,\dots,1)\in X:$ by the jacobian criterion of regularity we know that\quad $rank\ J_0(1) =k  $.\\
Further for $\ell\in \integ$, let $\ J_{\ell} \ $  denote  the submatrix    \  of 
$\ J_0(1)$ \ obtained by considering the rows\\ $\left(\displaystyle{\frac{\partial f_i }{\partial x_j}}(1,...,1)\right) $ \ with \ $d_i\in H_{\ell} $.
\end{remnot}   
    \begin{prop}\label{T11}
       \begin{enumerate}

 \item 
 {$ T_B^1=\bigoplus_{\ell\in\integsmall}T^1_B(\ell) $    is a $\integ$-$graded$
  $ \mathbb{ F}$-$vector\ space$ of finite dimension}:   
     \item[] ${\overline g}\in T^1_B(\ell)\II \textrm{deg} (g(f_j))= \textrm{deg} (f_j)+\ell\ \ \forall j$ \quad {\rm ({\it see} \cite{pi})}. 
\item
 { For $g\in T_B^1:$ \  \ $g=\displaystyle {\sum}_{i=1}^{k} \ \alpha_i t^{\beta_i}\displaystyle{\frac{\partial  }{\partial x_i}}$}.
 \item {$\textrm{dim}_  \mathbb{ F} T^1_B(\ell)=\n G_{\ell}-\textrm{dim}\ V_{\ell} -1$}, \  where \ $V_{\ell}$ is the sub-vector-space of $ \mathbb{ F}^{k+1} $ generated by the 
 row-vectors of the  matrix $J_{\ell}   $.
 \denu  
\end{prop}
  {\it Proof.} 2. We know that there exists $n\in \nat $ such that \  $   \m ^n T^1_B=0$.
  Therefore for each $g\in Hom_B(I/I^2,B)$   there exists $a\in \nat$ such that $t^{a}g\in Im\phi$. Further   in $B$  the Euler's identity holds:   
 $\displaystyle {\sum}_{i=0}^{k} \ n_i\ x_i\displaystyle{\frac{\partial f_j}{\partial x_i}}=0, \ \forall j=1,\dots, q  $  and so $g$ can be rewritten as a linear combination of the  partial derivatives  with respect to $(1,\dots,k)$.\\
3. Recall that $Im(\varphi)$ is generated by the the partial derivatives. For each    $ i\notin G_{\ell} ,$ we have:  $t^{ \ell+n_i }\in B$ and so  $t^{ \ell+n_i }\displaystyle{\frac{\partial\ }{\partial x_i}}\in Im\ \Phi $.  
On the contrary, note that \ for $i\in G_{\ell}$,  the vector  
$$\displaystyle{\sum}_{ 1}^k\alpha_i    t^{\ell+n_i}\displaystyle{\frac{\partial\ }{\partial x_i}}\in Hom_B(I/I^2,B)(\ell)   \II J_{\ell}\left( 0,\alpha_1, \alpha_2 ,\dots,\alpha_k  \right)^T=0, \ \  \forall \ d_j\in H_{\ell}.$$ 
Therefore the system has   $\infty^{\n G_l-\textrm{dim}\ V_{\ell}-1}$ solutions. 
\quad$\diamond$\\

For a semigroup $S$, let $S(1):=\{ n\in\integ\ |\ n+n_i\in S\ \ \forall i\geq 0\}$. 
  \begin{prop}\label{T10} Let $\ L=  S(1)\cup \{n\in\integ\ |\ n< -2c+2-2n_0\}$. Then \begin{enumerate}
  \item 
 $\textrm{dim}\  T^1_B(\ell)=0$ \ for each $\ell\in L$.
 \item in particular if $S$ is ordinary or hyperelliptic, then  $\textrm{dim}\  T^1_B(\ell)=0$ for each $\ell< -4g-2$.
  \item Let $f_1',\dots ,f_q'$ be a reordering of the set $\{f_1,\dots,f_q\}$ such that the degrees satisfy \\  $d'_1\leq d'_2\leq \dots \leq d_q'$. Let $J_0'(1) $ be the associated jacobian matrix   and let $p=minimum$ integer such that the first $p$ rows of $J_0'(1) $ constitue a matrix of rank $=k$. Then $T^1_B(\ell)=0$ for each $\ell<-d'_p$.
  \item $\textrm{dim}\  T^1_B(c-1-n_0-n_1)>0$ \quad $(see\ \cite{pi})$.
 \item If $\ell\geq c-2 n_1,$ then \ $\textrm{dim}\  T^1_B(\ell)=\textrm{max}\{0,\n G_{\ell}-1\}$
 \item If $\ell\geq c- 1-n_0,$ then \ $\textrm{dim}\  T^1_B(\ell)=0$.
 \item  For each $i=0,\dots,k$ \ we have: \ $ t^{c+n-n_i}\displaystyle{\frac{\partial  }{\partial x_i}}\in Im\ \Phi,\ \ \forall n\geq 0$.

 \denu
  \end{prop}
  {\it Proof.}  1. If $\ell\in S(1),$ then $G_{\ell}=\emptyset$ and we are done by (\ref{T11}.3)\\
  If $\ell<-2c+2-2n_0$, then $\n G_{\ell}=n_0$ and $n_i+n_j+\ell \leq 2(n_0+c-1)+\ell<0$, hence $H_{\ell}=\{1,\dots,q\} $. Then $\textrm{dim} \ V_{\ell}=n_0-1$ and the claim follows by (\ref{T11}.3).\\
  2. Follows by (1): if $S$ is ordinary, or hyperelliptic, then $-2c-2n_0 =-4g-2$.\\
  3. Immediate by the assumptions and by (\ref{T11}.3), since   we have: $\textrm{dim} \ V'_{\ell}=n_0-1$, because $H_\ell\supseteq\{1,\dots,p\}$.\\
   4.   Let $\ell=c-1-n_0-n_1:$ we have $\{ 0, 1\} \subseteq G_\ell$, \ while $H_{\ell}=\emptyset$ since $d_i>n_0+n_1 \ \ \forall i $.
Therefore the claim follows by (\ref{T11}.3). \\
5. Follows by (\ref{T11}.3) since in this case $H_{\ell}=\emptyset$.\\
6. Follows by (\ref{T11}.3) as a particular case.\\
7. Recall that $t^{c+n}\in B$ and $\textrm{deg}\left(\displaystyle{\frac{\partial  }{\partial x_i}}\right )=-n_i$.
\quad$\diamond$\\

 \noindent For  flatness conditions an essential  fact    is the following:
\begin{prop}\label{flat} $($See, e.g.   {\rm \cite[Page 8]{st}}$)$
  Given a cartesian diagram $$\begin{array}{cccc} \X&\hookrightarrow &Y\\
\Big\downarrow&&\ \ \Big\downarrow \pi\\
\{0\}&\hookrightarrow &\Sigma\end{array}$$ let $f_i  $ \  and   \ $F_i,\ i=1,\dots,q $  be respectively the equations of X  and $Y$.  Then:\\
the map $\pi $ is flat $\II$ every relation $\sum_1^q r_if_i=0,\ \  r_i,f_j\in  { \mathbb{ F}}[x_0,\dots,x_k]$ can be lifted to a relation $\sum_1^q R_iF_i=0, \ \ R_i,F_j\in  { \mathbb{ F}}[x_0,\dots,x_k]\otimes {\cal O}_{\Sigma,0}.$ \end{prop}

   \begin{thm} {\rm (See, e.g., \cite{st})}  The infinitesimal deformations  are in one-to-one correspondence with \ \  $Hom_B(I/I^2,\ B)$ as follows 
{$$\left(\begin{array}{ll}g:I/I^2\longrightarrow B\\
{  f}_j\mapsto {  g}_j(\textrm{mod}\ I)\\ \ \ (j=1,...,q)
\end{array}\right)\begin{array}{c}corresponds\\
to \ the\\
deformation\end{array}\ \ {  F}=\left(\begin{array}{c} f_1+\varepsilon g_1\\ ...\\ f_ q+\varepsilon g_ q\end{array} \right )$$} 
\end{thm}
{\it Proof.} (Outline) 
   The trivial deformations  (i.e. $Y\simeq \Sigma\times X$) lie in $Im\ \varphi$   \ $(\ref{T1})$.
   
  {\rm In fact these deformations   are   such that the ideal generated by the $(f_i+\varepsilon g_i)\in  { \mathbb{ F}} [\varepsilon, x_0,\dots,x_{k}]$ becomes equal, after a change of variables, to the ideal generated by the $(f_i)$.
 Now note that   a change of variables is $\{x_i\mapsto x_i+\varepsilon h_i\}$; since $\varepsilon^2=0$ easily one can see that it adds to each $g_i$ an element of the form $\ \sum_{j=0}^{k}\displaystyle{\frac{\partial f _i }{\partial x_{j}}}h_j$.\ \
 Therefore $T^1_B$ can be naturally identified with the set of infinitesimal deformations modulo the trivial ones.}\quad $\diamond$}\\
 
 Several semigroups have been recognized to be Weierstrass by  means of  the above theory: we collect in the following theorem the most important statements.
   \begin{thm} Assume  $X$ be an affine curve.
 \begin{enumerate}\item
 If $X$ is a complete intersection then $X$ is smoothable {\rm \cite{sc1}}. 
 \item If $X\subseteq \af^3$ or  $X$ is a Gorenstein curve of embedding dimension 4, then $X$   is smoothable $\cite{sm}, \ \cite{be}$. 
\item   Let \ $e,\  g$ denote respectively the multiplicity and the genus of the semigroup  $S$ and let $X= \textrm{Spec}( \mathbb{F}[S])$. Then:
\begin{enumerate}
\item     If \ $e\in\{3,4,5\}$, \ then $S$ is Weierstrass:
for $e=3$,  see also $\cite{mac}$, \   for $e=4$,  $e=5$, see ({\rm\cite{k1}, 
\cite {k2}}).\vspace{-0.1cm}
\item   If \ $g\leq 8$, then S is Weierstrass (\cite{k3}).

\item   If \   $2e>c-1$  and $g=9$, \ then S is Weierstrass  (\cite{k4}). 
\denu
\item   Let $H=\nat\setminus S$, define $\ \textrm{weight} (S):=\sum_{i=1}^g h_i-i$: 
if $\ \textrm{weight}  (S)\leq g/2$, \ then $S$ is Weierstrass ($ \cite{eh}$). 
\item    If \ $B$ is negatively graded $($i.e. $T^1_B(\ell)=0$ for each $\ell\geq 0)$, \ then $S$ is Weierstrass ($\cite{rv}$). 
\end{enumerate}
\end{thm}
 \subsection{Construction of the versal deformation with $\mathbb{G}_m$-action.}   
 With the above notations for the monomial curve $X:=\textrm{Spec}(\mathbb{ F}[S])$, \ $S$ a numerical semigroup,   we shall describe  Pinkham's {algorithm} \cite{pi} to construct  a  deformation $Y$ admitting a $\mathbb{G}_m$-action. Starting from the infinitesimal deformation associated to $\ \bigoplus _{\ell<0}T^1_B(\ell)$, by means of a finite number of steps one can obtain such deformation (with the greatest parameter space). Each step consists in the lifting of  a deformation on $\ \Sigma= \textrm{Spec} \  \mathbb{ F}[\varepsilon]/(\varepsilon)^n$ \ \   to a deformation on
$\ \Sigma'= \textrm{Spec} \  \mathbb{ F}[\varepsilon]/(\varepsilon)^{n+1}$.\\
Further in the last step we recall Pinkham's  construction (when possible) of a projective regular curve $\cal{C}$ admitting $S$ as semigroup at the point $P_{\infty}$ (see \cite[13.3]{pi}). This construction is  the main ingredient for the proof of Theorem \ref{fond}.\\

{\bf Step (0)} The first step of the algorithm is  the explicit computation of a $ \mathbb{ F}$-basis $E$ for \  $T_B^1$.\\

 {\bf Step (1) }
Let $ {r}$ be a $(p\times  q)$ matrix of relations among the generators $\{f_i\}$ of   $I$.  \\
For each $ {g}_j\in E$ construct a $(p\times  q)$ matrix  ${  \rho}_j={  \rho}_j(x_0,...,x_k)$,  such that ${  R}={  r}+\varepsilon {  \rho}_j$ is a relation matrix among the 
equations  of ${  F}={  f}+\varepsilon{  g}_j$, \ \ i.e.,
\\ \centerline{$({  r}+\varepsilon {  \rho}_j)({  f}+\varepsilon{  g}_j)={  r }{  f }+\varepsilon ({  r }{  g }_j+{  \rho }_j{  f })=0 \quad (\textrm{mod} \ \varepsilon^2).$}
A matrix $\rho_j$ such that  ${  \rho}_j{  f }=-{  r }{  g}_j $ exists since   any $g\in Hom_B(I/I^2,B)$ is a derivation (\ref{T1}.2), and so 
   the matrix $
   {r} {g}$\ has\ entries $\in I$,  \ for \ each \   $ {g}\in Hom_B(I/I^2,B).$
  In fact if $\sum  r_if_i=0, $ then   
 $0=g(\sum  r_if_i )=\sum r_ig(f_i)+\sum g(r_i) f_i =  \sum r_ig(f_i)\ \textrm{mod}\ I,\ \ i.e.\  \ \sum r_ig(f_i)\in I$.\\
 Hence
 any relation among  the $(f_i)$ lifts to a relation $R$ among the $(F_i)$, so that the projection $\pi$ is flat (\ref{flat}).
  \noindent Let $ E=<{  g}_1,...,{  g}_{m}>$ be the   $\mathbb{ F}$-basis of
 $\ \bigoplus _{\ell<0}T^1_B(\ell)$: assign a parameter $U_j$ to each $g_j$ with
    $$\textrm{weight}(U_j):= -\textrm{deg} (g_j). $$
     \  We obtain homogeneous    equations $$ {F}={  f}+\varepsilon({  g}_1U_1+...+{  g}_{m}U_{m})  \in   \mathbb{ F}[U_1,...,U_{m},x_0,...,x_k]$$  
for  a deformation $Y_1$ of $X$ with base space $\textrm{Spec}\  \mathbb{ F}[\varepsilon]/(\varepsilon)^2$.\\
By linearity the matrix { ${  \rho}:=U_1{  \rho}_1+...+U_{m}{  \rho}_{m}$ } \ \ is such that $   {r}+\varepsilon {  \rho}$ is a relation matrix for ${  F}$.\\

{\bf Step (2) } Now, called { ${  g}:=({  g}_1U_1+...+{  g}_{m}U_{m})$}, \ 
look for a vector ${  h}$ and for a matrix $ { \rho'}$ 
such that \\ \centerline{    $F={  f}+\varepsilon{  g}+\varepsilon^2{  h}$  \ \ and \ \    $  { {  R}}={  r}+\varepsilon{  \rho}+\varepsilon^2{  \rho'}$  }
verify 
\vspace{-0.1cm} 
 $${ {  R}}{ {  F}}={  r}{  f}+\varepsilon ({  r}{  g}+{ \rho}{  f})+\varepsilon^2({ \rho}{  g}+{  r}{  h}+ { \rho'}{  f})\equiv 0\quad \textrm{mod} \ (\varepsilon)^3 \vspace{-0.1cm}  $$ 
$${ {  R}}{ {  F}}=\varepsilon^2({ \rho}{  g}+{  r}{  h}+ { \rho'}{  f})\equiv 0\quad \textrm{mod} \ (\varepsilon)^3 \vspace{-0.1cm}  $$  
 Note that { ${ \rho}{  g}$ is quadratic in $U_1,...,U_{m}$}, therefore both ${ \rho'}, {  h}$ will be quadratic in $U_1,...,U_{m}$.
 To solve this equation we must impose several conditions to the variables $\{U_1,...,U_{m}\}$, but a solution exists since $X$ has a versal deformation by (\ref{pis}).\\  
 
   ....\\
  
 { \bf Step (n) } The matrices to find have entries of $degree  \ n$ in   $ U_1,...,U_{m}$. We already know that the algorithm ends.  Surely it ends when $\textrm{deg}(U_{i_1}...U_{i_h})> \textrm{deg} (f_j) \ \ \forall \ j $ and $
 \ \forall \ (i_1,...,i_h).$
 In fact at this step the needed matrices  are null by the theorem of existence of a   versal deformation for $X$ admitting a $ \mathbb{G}_m$-action \cite{pi}.\\
    
 { \bf Last Step  }  Let $R:=\mathbb{F}[U_1,...,U_n]/J$, \  $\Sigma= \textrm{Spec} (R )$ be the parameter space of the constructed deformation $Y$ of $X$ with $\mathbb{G}_m$-$action$ and let 
  $F=f+U_1{  g}_1+...+U_{m}{ g}_{m}+U_1^2h_{11}+\dots$ be the defining equations of $Y$. Substitute $U_i x_{k+1}^{\textrm{weight}(U_i)}$ for $U_i$  and let $A:=R[x_0,\dots,x_{k+1}]/(F)$. Then the morphism $\pi:Proj(A)\longrightarrow \Sigma$ is flat and proper with fibres reduced projective curves \cite[13.4]{pi}. The generic fibre $\cal {C}$,   has only one regular point $P_{\infty}(t^{n_0}:t^{n_1}:\dots:t^{n_k}:0)$ at infinity. If one fibre $\cal {C}$ is regular, then the semigroup associated to the pair $(\cal {C},P_{\infty})$ is clearly equal to the semigroup $S$. 
 \section{Examples}
 In this section we show   the above algorithm in some particular example.  
  
   \subsection{The  case of embedding dimension 3}
   First we calculate explicitely a deformation with $\mathbb{G}_m$-action for a monomial curve $X\subseteq \mathbb{A}^3_{\mathbb{F}}$.
    \begin{ex} \vspace{-0.1cm}
 Let $S=<4,9,11>$, \ $B=  \mathbb{ F}[S]$, \ $X:= \textrm{Spec}\ B$.  \\
The conductor is $c=15$, the  Apery set is\ \ ${\cal A}=\{n_0 =4,n_1 =9,n_ 2=11,n_3 =18\}$.  \\ 	The equations defining the curve $X$ in $\mathbb{F}[x_0,x_1,x_2]$ \ are $$f_1=x_0^5 - x_1x_2,\quad f_2= x_0x_1^2 - x_2^2,\quad f_3= -x_1^3 + x_0^4x_2 \vspace{-0.3cm}$$  
with matrix of relations: $ r= \left( \begin{array}{rrr }
	   -x_2& x_1& x_0\\   x_1^2&-x_0^4& -x_2\end{array}\right)=\left( \begin{array}{ccclllllllll}
	    r_1\\  {r}_2\end{array}\right)$ and  Jacobian matrix \\
	    $$  J_0= \Big(  
\ \    x_0    {\frac{\partial f_j }{\partial x_{0}}},x_1  {\frac{\partial f_j }{\partial x_{1}}}  \ ,x_2   {\frac{\partial f_j }{\partial x_{2}}} \Big)=    \left(\begin{array}{crrrcccccccccccccccccc}
  
  5x_ 0^{5}&-x_1x_2&-x_1x_2\ \\ 
  x_0x_1 ^{2} &2x_0x_1^2&-2x_2^2    \\
    4x_0^4x_2&-3x_1^3& x_0^4 x_2
   \end{array}  \right) $$
   $$  J_0\equiv \left (\begin{array}{lll}
      t^{20}&0&0    \\
   0&t^{22}&0  \\
 0&0&t^{27}  \\
  \end{array}  \right)\left(\begin{array}{rrrrcrrrrrccccccccccc}
   5&-1&-1\\
   1&2&-2\\
   4&-3&1
  \end{array}  \right)\qquad (\textrm{mod}\ I).$$ 
   \end{ex}

  Let { $\Delta _i:=x_i\frac{\partial  }{\partial x_{i}}, \ i=0,1,2$} (degree 0 derivations).\\
  
 {\bf Step (0)} One can easily see that $T^1(B)$ is generated as $B$-module by \\
  $$\begin{array}{ccc} T^1(-18):& t^{-18}(\Delta_1-\Delta_2):=D_1\\
  T^1(-16):& t^{-16}(\Delta_1+\Delta_2):=D_2\\
    T^1(-11): & t^{-11}(\Delta_1+\Delta_2):=D_3.
    \end{array}$$  with images the classes $\textrm{mod}\ I$ of 
  {$$<g_1=\left( \begin{array}{rrr}
0\   \\ 4x_0\\- 4x_1
\end{array}\right),\ g_2=\left( \begin{array}{rrr}
-2x_0  \\ 0\  \\ -2x_2
\end{array}\right),\ g_3=\left( \begin{array}{rrrllllllll}
-2x_1  \\ 0\  \\-2x_0^4
\end{array}\right)>$$}
( Note that  as $\mathbb{F}$-vector spaces we have
  $\textrm{dim}_ \mathbb{ F}T^1(B)=17,\ \ \textrm{dim}_ \mathbb{ F}T^1(B)^{-}=15$).\\
 {\bf Step (1) } Using the above algorithm (restricted to three generators) we get the infinitesimal   deformation\\ 
   \centerline{ $\pi: \textrm{Spec}\big(  \mathbb{ F}[\varepsilon]/(\varepsilon )^2\otimes   \mathbb{ F}[x_0,x_1,x_2]/I_1\big)\longrightarrow \textrm{Spec} \big(\mathbb{ F}[\varepsilon]/(\varepsilon )^2\Big)$,} 
 with $U_i\in  \mathbb{ F}$, \  $I_1$   generated    by  the rows of  $F_1={  f}+\varepsilon{  g},\ \  $ with \ \ ${  g}=   U_1{  g}_1+U_2  {  g}_2+U_3 {  g}_3$:
 $$ {  F_1}  = {  f}+\varepsilon\left[
 U_1\left( \begin{array}{ccc}
0  \\ x_0\\- x_1
\end{array}\right)+U_2 \left( \begin{array}{ccc}
x_0  \\ 0\\ x_2
\end{array}\right)+U_3 \left( \begin{array}{cccllllllll}
x_1  \\ 0\\ x_0^4
\end{array}\right)\right] \!\! ,\   \textrm{weight} (U_1,U_2,U_3)=(18,16,11).$$
 In fact there exists  the matrix  $\rho= \left( \begin{array}{rrrlllllllll}
	    -U_3   & 0 \ \  & 0\ \ 
	   \\  \   U_1  & -U_2  &  U_3  \end{array}\right)$ such that
  { $( r+\varepsilon \rho)(f+\varepsilon g)\equiv 0$} 
  $(\textrm{mod}\ (\varepsilon)^2$, i.e., $(rg+\rho f=0)$ (this assures      $\pi$ is flat, with $R_1:=r+\varepsilon \rho \ $   relation matrix for $F_1$):
 $$ rg= \left( \hspace{-0.2cm}\begin{array}{crrlllllllll}
	   U_3( x_0^5- x_1x_2)\\   
	  U_1(-x_0^5+x_1x_2)+U_2(x_0x_1^2-x_2^2)+U_3(-x_1^3+x_0^4x_2) \end{array} \right)=  - {\rho f} .$$
   {\bf Step(2)}  \   Now look for \ $h,\rho'$ \ such that { $F_2=f+\varepsilon g+\varepsilon^2 h$  \ and \  { $R_2=r+\varepsilon \rho+\varepsilon^2 \rho'$}}
	 satisfy \  $F_2R_2=0\ (\textrm{mod} \ (\varepsilon)^3)$, i.e., $\rho g+rh+ \rho'f\equiv 0$. Get     $$ {\rho g} = \left( \begin{array}{crrlllllllll}
	    -U_2U_3x_0-U_3^2x_1  	   \\  
	    -U_1U_2x_0-U_1U_3x_1+U_1U_2x_0+U_2U_3x_2+U_1U_3x_1+U_3^2x_0^4 \end{array}\right)=$$
	   	    $ = \left( \begin{array}{rrrlllllllll}
	   -x_2    & x_1\   & x_0\   
	   \\  \ x_1^2 &-x_0^4   & -x_2  \end{array}\right)
	   \left( \begin{array}{c}
	  0\\
		  -U_3^2 \\
		    -U_2U_3  \\
 \end{array}\right)=-rh,$  \ \
 with \ \ {$h=  \left( \begin{array}{c}
	  0\\
		   U_3^2 \\
		    U_2U_3  \\
 \end{array}\right)$}. \\
 Finally one can see that  \ {${\rho h}=0,$}  \ \ therefore we can choose $\rho'=0$.\quad Hence the algorithm ends at the second step and a deformation of $f$ on $\textrm{Spec}\  \mathbb{ F}[U_1,U_2,U_3]$  has homogeneous weighted equations 
 { $$ {  F}  = {  f}+  
 \left( \begin{array}{ccc}
U_2 x_0 +U_3x_1 \\ 
U_1x_0\\
- U_1x_1+U_2x_2+U_3x_0^4
\end{array}\right)+\left( \begin{array}{c}
	  0\\
		   U_3^2 \\
		    U_2U_3  \\
 \end{array}\right).$$}
\begin{rem} {\rm Note that in the entries of the matrix $h$ the coefficient of $U_1^2$ is null. This is clear since   $\ \textrm{deg}(U_1^2)=36> \textrm{deg} (f_i),\ \forall\ i=1,2,3$, \  and the  equations are homogeneous according to the existence of a  $ \mathbb{G}_m$-action. Hence if we restrict to $g_1$,   we get the    deformation    \\ 
   \centerline{$\pi: Y=  \textrm{Spec} \big (  \mathbb{ F}[U_1] \otimes   \mathbb{ F}[x_0,x_1,x_2]/J\big)\longrightarrow  \mathbb{ A}^1_ {\mathbb{ F}}$ }    
 with    \ the ideal \ $J$   generated    by  the rows of   $$ {  F_1}  =   \left( \begin{array}{ccc} x_0^5 - x_1x_2\\   x_0x_1^2 - x_2^2\\ -x_1^3 + x_0^4x_2\end{array}\right)+
 U_1\left( \begin{array}{ccc}
0  \\ x_0\\- x_1
\end{array}\right)=   \left( \begin{array}{ccc} x_0^5 - x_1x_2\\ 
   x_0x_1^2 - x_2^2+U_1x_0\\ -x_1^3 + x_0^4x_2-U_1x_1\end{array}\right) .$$ 
      The algorithm ends at step (1) with smooth parameter space $\mathbb{ A}^1_{\mathbb{ F}} $.
   The Jacobian matrix of the generic  fiber of $\ \pi\ $ is 
  $$  \left [\begin{array}{cccrcccccccccccccccccc}
  
  5x_ 0^{4}&- x_2&- x_1\ \\ 
   x_1 ^{2}+U_1 &2x_0x_1 &-2x_2     \\
    4x_0^3x_2&-3x_1^2-U_1& x_0^4 
   \end{array}  \right] $$
   One can check that the generic fiber is non singular.}
   \end{rem}

 {\bf The general 3-Space case.}  By means of a costruction due to  Patil-Singh \cite{psing} we can compute directly the equations of the monomial curve associated to a semigroup $S=<n_0,n_1,n_2>$. 
 In this case we already know that every semigroup singularity is smoothable by Shaps' paper \cite{sm}: here the equations of a deformation are obtained as minors of a suitable matrix.
Let $S=<n_0,n_1,n_2>$, \ $n_0<n_1<n_2$,   let $Ap(S)$ be the Apery set  respect to $n_0$ and let  
$$
\left\{\begin{array}{ll}u:=\textrm{min}\{n\in\nat\ |\ un_1\in<n_0,n_2>,\ \ un_1\notin Ap(S)\}\\ 
v:=\textrm{min}\{n\in\nat\ |\ vn_2\in<n_0,n_1> \}
\end{array}
\right.
$$
Then 
$\left\{\begin{array}{ll}un_1=\lambda n_0+wn_2,\ \lambda\geq 1\\
vn_2=\mu n_0+zn_1,\ \ v\geq 2,\ \ v>w,\ 0\leq z<u\\
\quad \textrm{further:}\\
(\lambda+\mu)n_0=(u-z)n_1+(v-w)n_2.
\end{array}\right.\qquad(*)$.\\
By  \cite{psing} we know that the curve is a complete intersection $\II zw\mu=0$.\\
Then assume $zw\mu\neq 0$: we get the following      generators for the ideal $I$ and the relation module $r$ of   $X$:
$$I=\left\{\begin{array}{lll}f_1=x_1 ^{u }-x_ 0^{\lambda }x_2 ^{w }\\
f_2=x_1 ^{u-z }x_2 ^{v-w }-x_0 ^{\lambda+\mu  }\\
f_3=x_2 ^{v }-x_0 ^{\mu }x_1 ^{ z}
\end{array}\right.;\quad r=
\left(\begin{array}{ccc}  -x_2 ^{v-w } & x_1 ^{z } &\  -x_0 ^{\lambda } \\
  x_0 ^{ \mu} & -x_2 ^{w } & \ x_1 ^{u-z } \end{array}\right) =\left(\begin{array}{ccc}  r_1\\
  r_2 \end{array}\right) .$$
  Let $e_i$ denote the $i-th$ unit row vector. By Shaps' algorithm we get the following set of generators  of $Hom(I/I^2,B)$ as a $B$-module:
$$h_{11}:\left\{\begin{array}{lll}
f_1\mapsto \textrm{det}(e_1,e_1,r_2)=&0\\
 f_2\mapsto \textrm{det}(e_1,e_2,r_2)=&x_1^{u-z}\\
 f_3\mapsto \textrm{det}(e_1,e_3,r_2)=&x_2^{w}
 \end{array}\right.
 \quad h_{12}:\left\{\begin{array}{lcl}
f_1\mapsto \textrm{det}(e_1,r_1,e_1)=&0\\
 f_2\mapsto \textrm{det}(e_1,r_1,e_2 )=&x_0^{\lambda}\\
 f_3\mapsto \textrm{det}(e_1,r_1,e_3 )=&x_1^z
 \end{array}\right.$$
 $$h_{21}:\left\{\begin{array}{lcl}
f_1\mapsto \textrm{det}(e_2,e_1,r_2)=&-x_1^{u-z}\\
 f_2\mapsto \textrm{det}(e_1,e_2,r_2)=&0\\
 f_3\mapsto \textrm{det}(e_1,e_3,r_2)=&x_0^{\mu}
 \end{array}\right.
 \quad h_{22}:\left\{\begin{array}{lcl}
f_1\mapsto \textrm{det}(e_2,r_1,e_1)=&-x_0^{\lambda}\\
 f_2\mapsto \textrm{det}(e_2, r_1,e_2 )=&0\\
 f_3\mapsto \textrm{det}(e_2,r_1,e_3 )=&x_2^{v-w}
 \end{array}\right.$$
 $$h_{31}:\left\{\begin{array}{lcl}
f_1\mapsto \textrm{det}(e_3,e_1,r_2)=&-x_2^{w}\\
 f_2\mapsto \textrm{det}(e_2,e_2,r_2)=&-x_0^{\mu}\\
 f_3\mapsto \textrm{det}(e_1,e_3,r_2)=&0
 \end{array}\right.
 \quad h_{32}:\left\{\begin{array}{lcl}
f_1\mapsto \textrm{det}(e_3,r_1,e_1)=&-x_1^z\\
 f_2\mapsto \textrm{det}(e_3, r_1,e_2 )=&-x_2^{v-w}\\
 f_3\mapsto \textrm{det}(e_3, r_1,e_3 )=& 0
 \end{array}\right.$$
 We can construct the infinitesimal deformation (not miniversal, since $\textrm{dim}\ T^1_B$ is   greater,
  in general, but the other generators as vector space have greater degrees ). $${F }={f }+\epsilon\left[U_1\left(\!\!\begin{array}{lll}
 0\\
  x_1^{u-z}\\
  x_2^{w}
 \end{array}\!\!\right)+U_2\left(\!\!\begin{array}{lcl}
 0\\
 x_0^{\lambda}\\
 x_1^z
 \end{array}\!\!\right)+U_3\left(\!\!\begin{array}{cl}
 -x_1^{u-z}\\
  0\\
 x_0^{\mu}
 \end{array}\!\!\right)
+U_4\left(\!\!\begin{array}{cl}
 -x_0^{\lambda}\ \\
 0\\
 \ \ x_2^{v-w}
 \end{array}\!\!\right)+\right.$$ 
 $\left.\ \hspace{2.2cm}+U_5\left(\!\!\begin{array}{cl}
 -x_2^{w}\\
 -x_0^{\mu}\\
 0
 \end{array}\!\!\right)
 +U_6\left(\!\!\begin{array}{cl}
 -x_1^z\\
 \ \ x_2^{v-w}\\
  0
 \end{array}\!\!\right )\right]={f }+\epsilon{g} 
 $.\\
 With:\qquad $\textrm{weight}(U_1,...,U_6)=((v-w)n_2,\mu n_0, zn_1,wn_2,\lambda n_0, (u-z)n_1)$.  \\
 A relation matrix for $F$ is ${R}= {r}+\epsilon {\rho}$, \ 
 with $\ {\rho}=\left(\begin{array}{rrr}
 -U_1&-U_3&-U_5\\
 -U_2&-U_4&-U_6
 \end{array}\right)$. \\ 
 In fact
 $${rg}   =U_1\left(\!\!\begin{array}{lll}
 f_1\\0  \end{array}\!\!\right)+
 U_2\left(\!\!\begin{array}{lcl}
 0\\
 f_1
 \end{array}\!\!\right)+U_3\left(\!\!\begin{array}{cl}
f_2\\
  0 \end{array}\!\!\right)
+U_4\left(\!\!\begin{array}{cl}
 0\\
 f_2  
 \end{array}\!\!\right)+  U_5\left(\!\!\begin{array}{cl}
f_3\\
 0
 \end{array}\!\!\right)
 +U_6\left(\!\!\begin{array}{cl}
    0\\
    f_3
 \end{array}\!\!\right ).$$
 Now   the equation $({r}+\epsilon {\rho})({f }+\epsilon{g}+\epsilon^2{h} )=\epsilon({r} {g}   +{\rho}{f } )+\epsilon^2({\rho g }+{rh})=0$\\  has the solution  \\ \centerline{$ {h}=
 \left(\begin{array}{ccc}
 U_3U_6-U_4U_5\\
 U_2U_5-U_1U_6\\
 U_1U_4-U_2U_3
 \end{array}
 \right).$ } \vspace{0.1cm}\\
 Further the entries of  ${h}$ are the  $2\times 2$ minors of the matrix ${\rho}$ so that ${\rho}{h}=0$: hence there are no obstructions (conditions on $\{U_i\}$ necessary to have flatness).\\
 The  lift to a deformation with parameter space $\textrm{Spec} (\mathbb{ F}[U_1,...,U_6])$ is 
$${F }={f }+  U_1\left(\!\!\begin{array}{lll}
 0\\
  x_1^{u-z}\\
  x_2^{w}
 \end{array}\!\!\right)+U_2\left(\!\!\begin{array}{lcl}
 0\\
 x_0^{\lambda}\\
 x_1^z
 \end{array}\!\!\right)+U_3\left(\!\!\begin{array}{cl}
 -x_1^{u-z}\\
  0\\
 x_0^{\mu}
 \end{array}\!\!\right)
+U_4\left(\!\!\begin{array}{cl}
 -x_0^{\lambda}\ \\
 0\\
 \ \ x_2^{v-w}
 \end{array}\!\!\right)+ $$ 
 $$  
 +U_5\left(\!\!\begin{array}{cl}
 -x_2^{w}\\
 -x_0^{\mu}\\
 0
 \end{array}\!\!\right)
 +U_6\left(\!\!\begin{array}{cl}
 -x_1^z\\
 \ \ x_2^{v-w}\\
  0
 \end{array}\!\!\right ) +  \left(\begin{array}{ccc}
 U_3U_6-U_4U_5\\
 U_2U_5-U_1U_6\\
 U_1U_4-U_2U_3
 \end{array}
 \right).$$
Since $X$  is smoothable \cite{sm}, we deduce in particular that $(0,0,0)$ is a regular point on the general fibre: hence   
$$
1\in\{u-z,z,v-w,w,\lambda,\mu\}.
$$ 

\subsection{The example of Buchweitz.}

We show what happens in the following case of a non-smoothable monomial curve.

\begin{ex}\label{buc2}{\rm 
This example due to Buchweitz \cite{bu}  shows the first known case of non-smoothable monomial curve (see \cite{bu}). We calculate explicitly the miniversal deformation. 
Let 
$$
S=< 13, 14, 15, 16, 17, 18, 20, 22, 23>
$$
  The ring $B= \mathbb{ F}[S]$ has 32 equations in $\mathbb{ F}[x_0,\dots,x_8]$ \ (found by means of {\cocoa} \cite{cocoa}):
  $$ \begin{array}{rrrrr}-x_1^2+x_0x_2&
 -x_2^2 + x_1x_3&
 -x_1x_2 + x_0x_3&
 -x_3^2 + x_2x_4\\
 -x_2x_3 + x_1x_4&
 -x_1x_3 + x_0x_4&
  -x_4^2 + x_3x_5&
 -x_3x_4 + x_2x_5\\
 -x_2x_4 + x_1x_5&
  -x_1x_4 + x_0x_5&
 -x_5^2 + x_3x_6&
 -x_4x_5 + x_2x_6\\
 -x_3x_5 + x_1x_6&
 -x_2x_5 + x_0x_6&
 x_0^2x_1 - x_6^2&
 -x_0^2x_3 + x_6x_7\\
 -x_6^2 + x_5x_7&
 -x_0^3 + x_4x_7&
 -x_5x_6 + x_3x_7&
 -x_4x_6 + x_2x_7\\
-x_3x_6 + x_1x_7&
-x_2x_6 + x_0x_7&
x_0^2x_5 - x_7^2&
 -x_0x_1x_5 + x_7x_8\\
 -x_0^2x_4 + x_6x_8&
 -x_0^2x_2 + x_5x_8 &
  -x_5x_7 + x_4x_8&
 -x_4x_7 + x_3x_8\\ 
 -x_3x_7 + x_2x_8&
-x_2x_7 + x_1x_8&
 -x_1x_7 + x_0x_8&
-x_0^2x_6 + x_8^2\ \ \end{array}  $$
 The   Jacobian\ matrix whose $rank_P$ is $8$  if $P\in X,\ \ P\neq (0,\dots,0)$
is the following: 

$   J= \left[\begin{array}{ccccccccccc}
 x_2 &  -2x_1 &  x_0 &  0 &  0 &  0 &  0 &  0 &  0 \\
    0 &  x_3 &  -2x_2 &  x_1 &  0 &  0 &  0 &  0 &  0 \\
    x_3 &  -x_2 &  -x_1 &  x_0 &  0 &  0 &  0 &  0 &  0 \\
    0 &  0 &  x_4 &  -2x_3 &  x_2 &  0 &  0 &  0 &  0 \\
    0 &  x_4 &  -x_3 &  -x_2 &  x_1 &  0 &  0 &  0 &  0 \\
    x_4 &  -x_3 &  0 &  -x_1 &  x_0 &  0 &  0 &  0 &  0 \\
    0 &  0 &  0 &  x_5 &  -2x_4 &  x_3 &  0 &  0 &  0 \\
    0 &  0 &  x_5 &  -x_4 &  -x_3 &  x_2 &  0 &  0 &  0 \\
    0 &  x_5 &  -x_4 &  0 &  -x_2 &  x_1 &  0 &  0 &  0 \\
    x_5 &  -x_4 &  0 &  0 &  -x_1 &  x_0 &  0 &  0 &  0 \\
    0 &  0 &  0 &  x_6 &  0 &  -2x_5 &  x_3 &  0 &  0 \\
    0 &  0 &  x_6 &  0 &  -x_5 &  -x_4 &  x_2 &  0 &  0 \\
    0 &  x_6 &  0 &  -x_5 &  0 &  -x_3 &  x_1 &  0 &  0 \\
    x_6 &  0 &  -x_5 &  0 &  0 &  -x_2 &  x_0 &  0 &  0 \\
    2x_0x_1 &  x_0^2 &  0 &  0 &  0 &  0 &  -2x_6 &  0 &  0\\
    -2x_0x_3 &  0 &  0 &  -x_0^2 &  0 &  0 &  x_7 &  x_6 &  0 \\
    0 &  0 &  0 &  0 &  0 &  x_7 &  -2x_6 &  x_5 &  0 \\
    -3x_0^2 &  0 &  0 &  0 &  x_7 &  0 &  0 &  x_4 &  0 \\
    0 &  0 &  0 &  x_7 &  0 &  -x_6 &  -x_5 &  x_3 &  0 \\
    0 &  0 &  x_7 &  0 &  -x_6 &  0 &  -x_4 &  x_2 &  0 \\
    0 &  x_7 &  0 &  -x_6 &  0 &  0 &  -x_3 &  x_1 &  0 \\
    x_7 &  0 &  -x_6 &  0 &  0 &  0 &  -x_2 &  x_0 &  0 \\
    2x_0x_5 &  0 &  0 &  0 &  0 &  x_0^2 &  0 &  -2x_7 &  0\\
    -x_1x_5 &  -x_0x_5 &  0 &  0 &  0 &  -x_0x_1 &  0 &  x_8 &  x_7 \\
 -2x_0x_4 &  0 &  0 &  0 &  -x_0^2 &  0 &  x_8 &  0 &  x_6 \\
 -2x_0x_2 &  0 &  -x_0^2 &  0 &  0 &  x_8 &  0 &  0 &  x_5 \\
    0 &  0 &  0 &  0 &  x_8 &  -x_7 &  0 &  -x_5 &  x_4 \\
    0 &  0 &  0 &  x_8 &  -x_7 &  0 &  0 &  -x_4 &  x_3 \\
    0 &  0 &  x_8 &  -x_7 &  0 &  0 &  0 &  -x_3 &  x_2 \\
    0 &  x_8 &  -x_7 &  0 &  0 &  0 &  0 &  -x_2 &  x_1 \\
    x_8 &  -x_7 &  0 &  0 &  0 &  0 &  0 &  -x_1 &  x_0 \\
    -2x_0x_6 &  0 &  0 &  0 &  0 &  0 &  -x_0^2 &  0 &  2x_8  
\end{array}\right]$

Now we summarize the computation of $\textrm{dim}\ T^1_B(\ell)$ by means of the formula \vspace{0.1cm}
\\ \centerline{$\textrm{dim} \ T^1_B(\ell)=\n G_{\ell}-1-\rho_{\ell}.$}\vspace{0.1cm}
It is useful to consider the Jacobian matrix evaluated in $P(1,\dots,1)$ with the rows ordered by degree:
here the first column shows the  weighted degrees of the   equations.

 $$\Big(\ \textrm{deg}\ |\ J(1)\ \Big) =
\begin{array}{r|rrrrrrrrrrrrrrrrrrr}

 & x_0 & x_1& x_2& x_3& x_4& x_5& x_6& x_7& x_8\\
\hline
\textrm{deg}& 13&14&15&16&17&18&20&22&23\\ 
\hline
\hline 
 28&1 & -2 & 1 & 0& 0& 0& 0& 0& 0\\
29& 1 & -1 & -1 & 1 & 0& 0& 0& 0& 0\\
  30& 0& 1 & -2 & 1 & 0& 0& 0& 0& 0\\
30&  1 & -1 & 0& -1 & 1 & 0& 0& 0& 0\\
  31&  0& 1 & -1 & -1 & 1 & 0& 0& 0& 0\\
   31& 1 & -1 & 0& 0& -1 & 1 & 0& 0& 0\\

32&  0& 0& 1 & -2 & 1 & 0& 0& 0& 0\\
 32& 0& 1 & -1 & 0& -1 & 1 & 0& 0& 0\\
33&  0& 0& 1 & -1 & -1 & 1 & 0& 0& 0\\
   33&  1 & 0& -1 & 0& 0& -1 & 1 & 0& 0\\

34&  0& 0& 0& 1 & -2 & 1 & 0& 0& 0\\
34&  0& 1 & 0& -1 & 0& -1 & 1 & 0& 0\\

35&  0& 0& 1 & 0& -1 & -1 & 1 & 0& 0\\
  
35&  1 & 0& -1 & 0& 0& 0& -1 & 1 & 0\\
   36&  0& 1 & 0& -1 & 0& 0& -1 & 1 & 0\\
36&  0& 0& 0& 1 & 0& -2 & 1 & 0& 0\\
   36&  1 & -1 & 0& 0& 0& 0& 0& -1 & 1 \\
   37&  0& 1 & -1 & 0& 0& 0& 0& -1 & 1 \\

   37&  0& 0& 1 & 0& -1 & 0& -1 & 1 & 0\\
   38&  0& 0& 1 & -1 & 0& 0& 0& -1 & 1 \\

    38&  0& 0& 0& 1 & 0& -1 & -1 & 1 & 0\\
   39& -3 & 0& 0& 0& 1 & 0& 0& 1 & 0\\
   39&  0& 0& 0& 1 & -1 & 0& 0& -1 & 1 \\

40&  0& 0& 0& 0& 0& 1 & -2 & 1 & 0\\
40&  2   & 1 & 0& 0& 0& 0& -2 & 0& 0\\
   40&  0& 0& 0& 0& 1 & -1 & 0& -1 & 1 \\
   41&  -2   & 0& -1 & 0& 0& 1 & 0& 0& 1 \\

42&  -2   & 0& 0& -1 & 0& 0& 1 & 1 & 0\\
43&  -2   & 0& 0& 0& -1 & 0& 1 & 0& 1 \\

44&  2   & 0& 0& 0& 0& 1 & 0& -2 & 0\\
45&  -1   & -1   & 0& 0& 0& -1   & 0& 1 & 1 \\
46&  -2   & 0& 0& 0& 0& 0& -1 & 0& 2 
\end{array}$$
The matrix associated to degree $0$ derivations $\textrm{mod}\ I$ is
$$J_0=\left(x_i\displaystyle{\frac{\partial F_j }{\partial x_{i}}}\right)=
\left[\begin{array}{llll}
  t^{28} & 0&\ \dots \ &0 \\
 0& t^{29}&\ \dots  \ &0  \\
 \vdots&&\ddots\\
0 &0 &\ \dots\ &t^{46} 
\end{array}\right] J(1).$$
Now we show that $\textrm{dim}_ {\bf F}\  T^1_B=21$; to find a basis for   $   T^1_B(\ell) $ we have to solve the homogeneous system  associated to the minor of $J(1)$ formed by the rows of $\textrm{weight} \in H_{\ell}$.\\
\newpage 
{\bf Step (0)}\ \ First for each $\ell\in \integ$ we describe the subsets $G_{ell},\ H_{ell}$.

$$  \begin{array}{rlcllc} \ell\ \ \ & \ \ \ G_{\ell}&\n G_{\ell} &\ \ \ H&\ \ \rho&\textrm{dim} \ T^1_{\ell}\\
 -23\ \ &\ \{0,..,7\}&8 &\{28,..,35,42,44\} &\rho=7 & 0 \\
  -22\ \ &\ \{0,..,6,8\} &8& \{28,..,34,41,43,46\} &\rho=7 & 0  \\
  -21\ \ &\ \{0,...,8\} &9& \{28,..,33,40,42,45,46\}&\rho=8  & 0  \\
  -20\ \ &\ \{0,...,5,7,8\} &8& \{28,...,32,39,41,44,45\}&\rho=8  & 0  \\
  -19\ \ &\ \{0,...,8\} &9& \{28,...,31,38,40,43,44 \}&\rho=8 & 0  \\
  -18\ \ &\ \{0,...,4,6,7,8\} &8& \{28,...,30,37,39,42,43\} &\rho=8 & 0  \\
  -17\ \ &\ \{0,...,3,5,...,8\} &8& \{28,29,36,38,41,42\}&\rho=7  & 0  \\
  -16\ \ &\ \{0,1,2,5,6,7,8\} &7& \{28,35,37,40,41\}&\rho=7  & 0  \\
  -15\ \ &\ \{0,1,3,...,8\} &8& \{34,36,39,40\} &\rho=7 & 0  \\
  -14\ \ &\ \{0, 2,,...,8\} &8& \{33,35,38,39\}&\rho=7  & 0  \\
  -13\ \ &\ \{ 1,...,8\} &8& \{32,34,37,38\}&\rho=7 & 0  \\
  -12\ \ &\ \{0,...,8\} &9& \{31,33,36,37\}&\rho=7  & 1  \\
  -11\ \ &\ \{0,...,8\} &9& \{30,32,35,36 \}&\rho=7  & 1  \\
  -10\ \ &\ \{0,...,7\} &8& \{29,31,34,35\}&\rho=6  & 1  \\
  -9\ \ &\ \{0,...,6\} &7& \{28,30,33,34\}&\rho=5  & 1  \\
  -8\ \ &\ \{0,...,6\} &7& \{29,32,33\}&\rho=5  & 1  \\
  -7\ \ &\ \{0,...,5\} &6& \{28,31,32 \}&\rho=4  & 1  \\
  -6\ \ &\ \{0,...,5\} &6& \{30,31\}&\rho=4  & 1 \\
  -5\ \ &\ \{0,...,4\} &5& \{29,30 \}&\rho=3  & 1 \\
  -4\ \ &\ \{0,1,2,3,8\} &5& \{28,29\} &\rho=2  & 2 \\
  -3\ \ &\ \{0,1,2,7\} &4& \{28\} &\rho=1  & 2 \\
  -2\ \ &\ \{0,1,8\} &3& \emptyset& & 2 \\
  -1\ \ &\ \{0,6,7\} &3& \emptyset& & 2 \\
  1\ \ &\ \{5,6\} &2& \emptyset && 1 \\
  2\ \ &\ \{4\} &1& \emptyset && 0 \\
  3\ \ &\ \{3,5 \} &2& \emptyset && 1 \\
  4\ \ &\ \{2,4\} &2& \emptyset && 1 \\
  5\ \ &\ \{1,3\} &2& \emptyset& & 1 \\
  6\ \ &\ \{0,2\} &2& \emptyset& & 1  
  \end{array}$$
{\bf Step (1)}\ \  By using \lq\lq FreeMat" (see \cite{freemat}) we can construct the miniversal deformation \ (we present in detail the case $\ell=-12$ with $H_{\ell}=\{31,33,36,37\}$.\\Let $a$ be the submatrix of $J(1)$ formed by the rows with degrees $\in H_{\ell}$:\\ \ \ 
$$a=\begin{array}{r||rrrrrrrrrrrrrrrrrrr}
& x_0 & x_1& x_2& x_3& x_4& x_5& x_6& x_7& x_8\\
\textrm{deg}& 13&14&15&16&17&18&20&22&23\\
\hline
\hline
   31&  0& 1 & -1 & -1 & 1 & 0& 0& 0& 0\\
   31& 1 & -1 & 0& 0& -1 & 1 & 0& 0& 0\\

 33&  0& 0& 1 & -1 & -1 & 1 & 0& 0& 0\\
   33&  1 & 0& -1 & 0& 0& -1 & 1 & 0& 0\\

   36&  0& 1 & 0& -1 & 0& 0& -1 & 1 & 0\\
36&  0& 0& 0& 1 & 0& -2 & 1 & 0& 0\\
   36&  1 & -1 & 0& 0& 0& 0& 0& -1 & 1 \\
   37&  0& 1 & -1 & 0& 0& 0& 0& -1 & 1 \\

   37&  0& 0& 1 & 0& -1 & 0& -1 & 1 & 0\\
\end{array}$$  \newpage
\noindent {\bf Step (1.1)} \ \ 
\\ Write the matrix $b$  obtained by deleting the deg-column in $a$,   find $rank(b)$ and   a total reduction   $c$ of $b$, that is 
 $$
 c=\left[\begin{array}{rrrrrrrrrrrrr}
 
   1 &  0 &  0&   0  & 0 &  0 &  0 &-10 &  9\\
			
   0 &  1  & 0&   0&   0&   0&   0&  -9 &  8\\
 			 
  0   &0 &  1&   0&   0&   0&   0&  -8&   7\\
 			   0 &  0&   0&   1&   0&   0&   0&  -7&   6\\
 			 
  0  & 0&   0&   0&   1&   0&   0&  -6&   5\\
 		 
  0  & 0&   0&   0&   0&   1&   0&  -5&   4\\
  0  & 0&   0&   0&   0&   0&   1 & -3 &  2\\
  			 
 0&   0&   0&   0&   0&   0&   0&   0&   0\\
  
 0&   0&   0&   0&   0&   0&   0&   0&   0
\end{array}\right]$$ 
 {\bf Step (1.2)} \ \   Let { $\Delta _i:=x_i\frac{\partial  }{\partial x_{i}}, \ i=0,\dots, 8$} (degree 0 derivations). Find the degree-0 derivation whose coefficients are   a solution of the homogeneous system associated to $c$ and by using the Euler's identity, we obtain  a solution  where the  coefficient of $\Delta_0$ is null. Then: \vspace{0.2cm}\\
  $T^1_B(-12)=<t^{-12} \Delta(1)>,$ with   $\Delta(1):= 	\Delta_1+2\Delta_2+3\Delta_3+4\Delta_4+5\Delta_5+7\Delta_6+9\Delta_7+10\Delta_8.$ \vspace{0.2cm}\\
{\bf Step (1.3)} \ \   Let $  e:=[0,1,2,3,4,5,7,9,10]^T$, to obtain the image $g_1$ of $\Delta(1)$  make the product: 
\\
$ J(1) e=\left[\begin{array}{rrrrrrrrrrrrrrrrrrrrrr}   0,
   0,
   ...,
   0,
  13,
   0,
   0,
 -13,
   0,
  13,
  13,
  13,
 -13,
  13,
  13\end{array}\right]^T  $ and so    $\Delta(1)$ takes $f$ to \ 
 $ g_1= $ \\ $\left[  
  t^{28 },t^{29},t^{30},t^{30},t^{31},t^{31},t^{32},t^{32},
  t^{33},t^{33},t^{34},t^{34},t^{35},t^{35}, t^{36},t^{36},  t^{36},t^{37},
  t^{37},t^{38},t^{38},t^{39},t^{39},t^{40 },\right.$\\ $\left. t^ {40},t^{40 },t^{41 },t^{42 },
  t^{43 },t^{44},t^{45} ,t^{46}
  \right]*J(1)e=$\\ $
  \left[\begin{array}{rrrrcccccccc}   0,
   0,
   ...   0,
  13t^{27},
   0,
   0,
 -13t^{28},
   0,
  13t^{29},
  13t^{30 },
  13t^{31 },
 -13t^{32 },
  13t^{33 },
  13t^{34}\end{array}\right]^T=$\\
  $  \left[\begin{array}{ccccccc}   0,
   0,
   ...
   0,
  13x_0x_1,
   0,
   0,
 -13x_1^2,
   0,
  13x_1x_2,
  13x_2^2,
  13x_1x_4,
 -13x_1x_5,
  13x_2x_5,
  13x_4^2\end{array}\right]^T\in \Big(M^2\Big)^{32}$.\\
(Here $*$ denotes the $pairwise$ vector product). Analogously we have:\\ 
$T^1_B(-11)=t^{-11}<\Delta(1)>.$\\
 $f\mapsto g_2= 
  \left[\begin{array}{rrrrcccccccc}   0,
   0,
   ...   0,
  13t^{28},
   0,
   0,
 -13t^{29},
   0,
  13t^{30},
  13t^{31 },
  13t^{32 },
 -13t^{33 },
  13t^{34 },
  13t^{35}\end{array}\right]^T=$\\
  $  \left[\begin{array}{ccccccc}   0,
   0,
   ...
   0,
  13 x_1^2,
   0,
   0,
 -13x_1x_2,
   0,
  13x_2^2,
  13x_1x_4,
  13x_1x_5,
 -13x_2x_5,
  13x_4^2,
  13x_4x_5\end{array}\right]^T\in \Big(M^2\Big)^{32}$.
\\ $T^1_B(-10)=t^{-10} \Delta(2)\  with  \  \Delta(2) := \Delta_1+2\Delta_2+3\Delta_3+4\Delta_4+5\Delta_5+7\Delta_6+9\Delta_7 .$

With image $g_3=t^{-10} [ 
  t^{28 },t^{29},t^{30},t^{30},t^{31},t^{31},t^{32},t^{32},
  t^{33},t^{33},t^{34},t^{34},t^{35},t^{35}, t^{36},t^{36},  t^{36},t^{37},\\
  t^{37},t^{38},t^{38},t^{39},t^{39},t^{40 },   t^ {40},t^{40 },t^{41 },t^{42 },
  t^{43 },t^{44},t^{45} ,t^{46}] * 
[     0,   0,   0,   0,   0,   0,   0,   0,   0,   0,    0,   0,   0,   0,   0,   0,\\ 
  $-$10,  $-$10,   0, $-$10,   0,  13,$-$10,   0, $-$13, $-$10,   3,  13,   3, $-$13,   3,    $-$\left. 7 \right.]=[0,\dots,0, $-$10t^{26},$-$10t^{27},0 ,\dots ]$.\\
    Therefore ${  f} \mapsto {  g_3} \in \Big(M^2\Big)^{32}$.\quad In the same way one obtains: \\
 $T^1_B(-9)=t^{-10}<\Delta(3)>,$ \ with \ $\Delta(3)= 	\Delta_1+2\Delta_2+3\Delta_3+4\Delta_4+5\Delta_5+7\Delta_6   ;$ \\
  $T^1_B(-8)=t^{-8}<\Delta(3)>.  $  Further one can  find that\\
  $T^1_B(-7)=t^{-7}<\Delta_1+2\Delta_2+3\Delta_3+4\Delta_4+5\Delta_5>$
 \\    $T^1_B(-6)=t^{-6}<\Delta_1+2\Delta_2+3\Delta_3+4\Delta_4+5\Delta_5>$\\  $T^1_B(-5)=t^{-5}<\Delta_1+2\Delta_2+3\Delta_3+4\Delta_4 >$ \\
$T^1_B(-4)=t^{-4}<\Delta_1+2\Delta_2+3\Delta_3,\  \Delta_8 >$\\ 
$T^1_B(-3)=t^{-3}<\Delta_1+2\Delta_2 ,\  \Delta_7 >:$ 
 $T^1_B(-2)=t^{-2}<\Delta_1 ,\  \Delta_8 > $ \\
  $  T^1_B(-1)=t^{-1}<\Delta_6 ,\  \Delta_7 > $ \\
   $\ T^1_B(1)=t<\Delta_5   > $ \\
     $\ T^1_B(3)=t^3<\Delta_5   > $ \\
       $\ T^1_B(4)=t^4<\Delta_4   >$ \\
         $\ T^1_B(5)=t^5<\Delta_3  > $ \\
           $\ T^1_B(6)=t^6<\Delta_2   > .$ 

We conclude    that each generator of $ T^1_B $ sends $f\mapsto g  \in \Big(M^2\Big)^{32}. $ \\
Hence all  the hypersurfaces defined by the equations $F_i\in\mathbb{F}[x_0,\dots,x_8,\alpha_1,\dots,\alpha_{21}], \ i=1\dots,32$ of the miniversal deformation are singular at $P(0,0,\dots,0,\alpha_1,\dots,\alpha_{21} ) $. In particular every fibre $X_{\underline T  }$ is singular: this means that $X $ is non-smoothable.}
\end{ex} 

 \section{Arithmetic sequences of embedding dimension four.}
 In this section using the above algorithms we   prove that semigroups of embedding dimension four generated by an arithmetic sequence are Weierstrass. 
 
  First we recall how to find the generators for the ideal $I$ of the monomial curve associated to the semigroup. We refer to the paper \cite{ps} and we use the same notations. 
 \begin{notat}\label{arit}
{\rm Assume $S=<n_0,\dots ,n_p,n_{p+1}>$, with $n_i=n_0+id$ (minimal system of generators), and denote by   $Ap(S)$   the Apery set  respect to $n_0$.\\
 Let $a,b \in \nat$ such that $n_0=a(p+1)+b$, with $a\geq 1$, \ $0\leq b\leq p$.\\
 For each $t\in\nat,$ let $q_t,\ r_t$ with $1\leq r_t\leq p$ such that $t=q_tp+r_t$,  let
 $g_t:= q_tn_p+n_{r_t}$ 
  and let  \\ \centerline{$\left\{\begin{array}{ll}
  u:=\textrm{min}\{t\in\nat\ |\ g_t \notin Ap(S)\}\\
   v:=\textrm{min}\{n\in\nat\ |\ vn_{p+1}\in<n_0,...,n_p> \}\end{array}\right.$}
Then 
$\left\{\begin{array}{ll}
\ \ g_u=\lambda n_0+wn_{p+1},\ \lambda\geq 1\\
vn_{p+1}=\mu n_0+g_z,\ \ v\geq 2,\ \ v>w,\ 0\leq z<u\\
\quad \textrm{further} :\\
(\lambda+\mu)n_0=g_{u-z} +(v-w)n_{p+1}.
\end{array}\right.\qquad(*)$ \\

\noindent It is easy to see that $u=p+1,\ \lambda=w=1$ \ and \\
if $b=0$, \ then $\ z=0$.\\
If $\ b\geq 1$ :
$v=a+1,\ \ \mu=a+d\geq 2,\ \ z=p+1-b$\\ and  a minimal set of generators for the ideal $I$ is the union of the following sets: \vspace{0.2cm}\\
$\xi_{ij}=\left[\begin{array}{llllll}
 x_ix_j-x_0x_{i+j} & \textrm{if} & i+j\leq p,&1\leq i\leq j \\
 x_ix_j-x_{i+j-p}x_p & \textrm{if} & i+j> p, &1\leq i\leq j \leq p-1
 \end{array}\right.$\vspace{0.2cm}\\
 $\phi_i=x_{1+i}x_p-x_ix_{p+1}\ \   \textrm{with}  \ \  0\leq i\leq p-1$\vspace{0.2cm}\\
 $\psi_j= x_{b+j}x_{p+1}^{v-1}-x_0^{\mu}x_j  \ \ \textrm{with}  \ \ 0\leq j\leq p-b$\vspace{0.2cm}\\
 $\theta =x_{p+1}^v-x_0^{\mu}x_{p+1-b}.$\vspace{0.2cm}}
 \end{notat}
 Now we deal with the case $\ p=2$ ($\textrm{embdim} (S)=4$): here\vspace{0.2cm}\\
 $\{\xi_{ij}\}=\{\xi_{11}\}=\{x_1^2-x_0x_2\}, \ \ \textrm{with} \ \ \textrm{deg}(\xi_{11})=2n_1$,\vspace{0.2cm}\\
$\{\phi_{i}\}=\{ \phi_0,\phi_1\}=\{x_1x_2-x_0x_3,\ x_2^2-x_1x_3\},\ \ \textrm{deg}(\phi_0)=n_1+n_2,\ \textrm{deg}(\phi_1)=2n_2$,\vspace{0.2cm}\\
$\{\psi_{j}\}=\left[\begin{array}{lll}\{ \psi_0,\psi_1\}=\{x_1x_3^{v-1}-x_0^{1+\mu},\ \ x_2 x_3^{v-1}-x_0^{\mu}x_1\}&  \textrm{if} & b=1, \ \textrm{deg}(\psi_1)=\mu n_0+n_1 ,\\
\{ \psi_0 \}=\{x_2x_3^{v-1}-x_0^{1+\mu} \}&  \textrm{if}  & b=2,  \ \textrm{deg}(\psi_0)=(1+\mu)n_0,
\end{array}\right.$\vspace{0.2cm}\\
 $\theta =\left[\begin{array}{lll}
 x_{3}^v-x_0^{\mu}x_{p+1-b},& \textrm{deg}(\theta)=\mu n_0+n_{p+1-b}&  \textrm{if} \ \ b=1,2\\
 x_{3}^v-x_0^{\mu} ,& \textrm{deg}(\theta)=\mu n_0 &  \textrm{if} \ \ b=0\\
 \end{array}\right.$.\vspace{0.2cm}\\
 Hence the equations for the associated monomial curve in $\mathbb{A}^4$ are :
 $$(b=0)\left(\begin{array}{lllll}
  x_1^2-x_0x_2 \\
 x_1x_2-x_0x_3\\
  x_2^2-x_1x_3\\
   x_3^{v}-x_0^{\mu} \end{array}\right); (b=1) \left(\begin{array}{lllll}
  x_1^2-x_0x_2 \\
 x_1x_2-x_0x_3\\
  x_2^2-x_1x_3\\
   x_1x_3^{v-1}-x_0^{1+\mu} \\
   x_2 x_3^{v-1}-x_0^{\mu}x_1 \\
   x_{3}^v-x_0^{\mu}x_{2} \end{array}\right); 
     (b=2) 
  \left(\begin{array}{lllll}
  x_1^2-x_0x_2 \\
 x_1x_2-x_0x_3\\
  x_2^2-x_1x_3\\
   x_2x_3^{v-1}-x_0^{1+\mu} \\
   
   x_{3}^v-x_0^{\mu}x_{1}\end{array}\right).$$

 \begin{lemma}\label{as4} Assume $S=<n_0,n_1,n_2,n_3>,$ minimally generated by an arithmetic sequence. With notation fixed in \ (\ref{arit}) \ we have: 
 \begin{enumerate}
 \item $  T^1_B(-\mu n_0)=<t^{-\mu n_0}(\Delta_1+2\Delta_2+3\Delta_3)>$.  \item  Further in case $b= 2  $, we have
 \item[]\ $  T^1_B(-( v-1)n_3 )=<t^{-(v-1)n_3}(\Delta_1+2\Delta_2+3\Delta_3)>$ 
 \item[]\ $  T^1_B(-n_2 )=<t^{-n_2}(2v\Delta_1+(v+1)\Delta_2+2\Delta_3)>$ 
 
 \denu
 \end{lemma}
{\it Proof.} (1). \ \  With notations
 (\ref{GH}) and (\ref{arit}),assume \  $\ell=-\mu n_0$.  Further recall that $\mu\geq 2$. Hence $\n G_{\ell}=4.$  Now  proceed  separately according that $b=0,1,2$.\\

  ({\bf Case} $b=0$). Easily one can see that $\quad  \ H_{\ell}=\{2n_1, n_1+n_2,2n_2\}$.\\
  The degree $0$ \ Jacobian matrix in this case is   
  $$  
J(0)=\Big(\displaystyle{x_i\frac{\partial f_j }{\partial x_i}}\Big)=
\left(\begin{array}{ccccccccc}
-x_0x_2&2x_1^2&-x_0x_2&0\\
-x_0x_3&x_1x_2&x_1x_2&-x_0x_3\\
0&-x_1x_3&2x_2^2&-x_1x_3\\
-\mu x_0^{\mu}&0&0&vx_3^{v}\\
\end{array}\right).
$$
The evaluation of this matrix in $P(1,\dots,1)\in X$ is 
$$   J(1)=\left(\begin{array}{ccccccccc}
- 1 &2 &- 1\ \ &0\\
- 1 &1   &1   &- 1 \\
\ \ 0&-1\    &2 &-1   \\
- \mu   &0&0&v \\
\end{array}\right) $$
Then $\textrm{dim}(V_{\ell})=2$ and $\textrm{dim} (T^1_B(-\mu n_0))=4-2-1=1$. A vector $(0,a,b,c)^T$ such that $J(1)(0,a,b,c)^T$ has the first three entries null is $(0,1,2,3)^T$. We obtain that a basis of 
$$T^1_B(-\mu n_0)\ \ {\rm is} \ \ t^{-\mu n_0}(\Delta_1+2\Delta_2+3\Delta_3).$$

( {\bf Case} $b=1$). We have   
 $\ H_{\ell}=\left[\begin{array}{ll}\{2n_1, n_1+n_2\}&if\ v=2\\
\{2n_1, n_1+n_2,2n_2\}&if\ v>2\\
\end{array}\right..$  \quad In fact \\
$2n_1-\mu n_0=n_2+(1-\mu)n_0\notin S$ since $\mu\geq 2$ and   $\{n_i\}$ is a minimal set of generators,\\
$n_1+n_2-\mu n_0=3d-(\mu-2)n_0\notin S$, since $n_0+3d=n_3 $,  
$$2n_2-\mu n_0=3n_2-vn_3=(3-v)n_0+(6-3v)d=
\left[\begin{array}{rl}n_0,&if\ v=2\\
<n_0,&if\ v>2
\end{array}\right.$$
for any   other generator $f_j$ of the ideal $I$, obviously $\textrm{deg}(f_j)-\mu n_0\in S$.\\
The $ \textrm{deg}\ 0$ \ Jacobian matrix  is $$J(0)=\Big(\displaystyle{x_i\frac{\partial f_j }{\partial x_i}}\Big)=\left(\begin{array}{ccccccccc}
-x_0x_2&2x_1^2&-x_0x_2&0\\
-x_0x_3&x_1x_2&x_1x_2&-x_0x_3\\
0&-x_1x_3&2x_2^2&-x_1x_3\\
-(1+\mu)x_0^{1+\mu}&x_1x_3^{v-1}&0&(v-1)x_1x_3^{v-1}\\
-\mu x_0^{\mu}x_1&-x_0^{\mu}x_1&x_2x_3^{v-1}&(v-1) x_2x_3^{v-1}\\
-\mu x_0^{\mu}x_2&0&-x_0^{\mu}x_2&vx_3^{v } \\
\end{array}\right)$$
The evaluation of this matrix in $P(1,\dots,1)\in X$ is 
$$   J(1)=\left(\begin{array}{ccccccccc}
- 1 &2 &- 1\ \ &0\\
- 1 &1   &1   &- 1 \\
\ \ 0&-1\    &2 &-1   \\
-(1+\mu)  &1&0&(v-1) \\
-\mu   &-  1\ \ &1&(v-1)  \\
-\mu    &0&-  1\ \ &v \\
\end{array}\right) $$
In both cases we see that $\textrm{dim}(V_{\ell})=2$. Then $\textrm{dim} (T^1_B(-\mu n_0))=4-2-1=1$ and analogously to case $b=0$, we recover the same   basis for 
$T^1_B(-\mu n_0)$.\\

( {\bf Case} $b=2$)   \ As above: \ $H_{\ell}= \{2n_1, n_1+n_2,2n_2\} $ 
 because \\ \centerline{$2n_2-\mu n_0=n_1+2n_2 -vn_3=(3-v)n_0+(5-3v)d\notin S$  (it is $<n_0$).}\\

The degree $0$ Jacobian matrix in this case   is 
$$
J(0)=\left(\begin{array}{ccccccccc}
-x_0x_2&2x_1^2&-x_0x_2&0\\
-x_0x_3&x_1x_2&x_1x_2&-x_0x_3\\
0&-x_1x_3&2x_2^2&-x_1x_3\\
-(1+\mu)x_0^{1+\mu}&0&x_2x_3^{v-1}&(v-1)x_2x_3^{v-1}\\
-\mu x_0^{\mu}x_1&-x_0^{\mu}x_1&0&vx_3^{v } \\
\end{array}\right)
$$
this  matrix evaluated in $P(1,\dots,1)$ is 
$$ 
 J(1)=\left(\begin{array}{ccccccccc}
- 1 &2 &- 1 &0\\
- 1 &1   &1   &- 1 \\
0&-1   &2 &-1   \\
-(1+\mu)  &0& 1&(v-1) \\
-\mu  &-  1&0&v  \\
\end{array}\right)\hspace{2cm}
$$
Therefore $\textrm{dim}(V_\ell)=1$, a basis for 
$T^1_B(-\mu n_0)\  \ {\rm is}\   t^{-\mu n_0}(\Delta_1+2\Delta_2+3\Delta_3).$\\

\noindent(3). Let $\ell=-(v-1)n_3$. Then:  
$$G_{\ell}=
\left[\begin{array}{ll}
\{0,1,2 \}, &if \ v=2 \\
\{0,1,2,3\}, &if \ v>2 \\
\end{array}\right.,\ \ 
H_{\ell}=\left[\begin{array}{ll}\{2n_1 \}&if\ v=2\\
\{2n_1, n_1+n_2,2n_2\}&if\ v>2\\
\end{array}\right..$$
In fact: \ $2n_1-(v-1)n_3\leq 2n_1-n_3=n_0-d\notin S$\\
$n_1+n_2-(v-1)n_3= n_0\in S$, if $v=2$, $n_1+n_2-(v-1)n_3<0 $, if $v>2$.\\
$2n_2-(v-1)n_3= n_1\in S$, if $v=2$, $n_1+n_3-(v-1)n_3<0 $, if $v>2$.\\
for any   other generator $f_j$ of the ideal $I$ , obviously $\textrm{deg}(f_j)-\mu n_0\in S$.\\
In both cases we conclude that $\textrm{dim}(T^1_B(\ell))=1$, with basis 
$$\ t^{-(v-1) n_3}(\Delta_1+2\Delta_2+3\Delta_3).$$
Let now    $\ell=-n_2$. Then:  \quad $G_{\ell}=
\{0,1,3 \},\ \ 
H_{\ell}=\{vn_3\}.$\\
In fact assume $vn_3-n_2\in S$, i.e.,  $vn_3=\alpha n_0+\beta n_1+\gamma n_2+\delta n_3,$ with $\gamma\geq 1$, then $\delta=0$ by the minimality of $v$; 
 $\beta\geq 1\I $(since $vn_3= \mu n_0+n_1)  \quad \mu n_0+n_1=\alpha n_0+(\beta-1)n_1+(\gamma-1)n_2+n_0+n_3\I 
 (v-1)n_3\in<n_0,n_1,n_2>,$ contradiction.
 Then $\beta=0$ and so $vn_3=\alpha n_0+\gamma n_2$, impossible since the residues mod $n_0$ cannot be equal.
 We conclude that $\textrm{dim}(T^1_B(\ell)=1$ and a basis is 
$t^{- n_2}(2v\Delta_1+(v+1)\Delta_2+2\Delta_3)$.\\
 
In next theorem we prove that any semigroup generated by an arithmetic sequence with embedding dimension 4 is Weierstrass. Further we find  the equations   of a 1-parameter flat family  of smooth projective with only one point $P_{\infty}$ at infinity and   the semigroup associated  at $P_{\infty}$ equal  to $S$. This is done by using Pinkham's algorithm \cite{pi}.

\begin{thm}    With notation \ref{as4}, assume the semigroup $S=<n_0,n_1,n_2,n_3>$ minimally generated by an arithmetic sequence and let $X:=\textrm{Spec} (\mathbb{F}[S])$: then $X$ is smoothable and S is Weierstrass. More precisely there exists one deformation $Y$ of $X$ with smooth   generic fibres (projective curves) and parameter space $\mathbb{A}^1_\mathbb{ F}:$
 \begin{enumerate}
 \item \ If \  $b=0$, \ the equations  \  $F=f+Ux_4^{\mu n_0} \left(\begin{array}{cc} 0\\ 0\\ 0\\ 1  \end{array}\right) =\left(\begin{array}{lllll}
  x_1^2-x_0x_2 \\
 x_1x_2-x_0x_3\\
  x_2^2-x_1x_3\\
   x_3^{v}-x_0^{\mu}+Ux_4^{\mu n_0}  \end{array}\right)$\\
   
  define the required deformation $\pi:Y\longrightarrow \mathbb{A}^1_\mathbb{ F}$.  \item \ If \  $b=1$, \ the equations  \  $F=f+Ux_4^{\mu n_0} \left(\begin{array}{cc} 0\\ 0\\ 0\\ x_0\\ x_1\\x_2 \end{array}\right) =\left(\begin{array}{lllll}
  x_1^2-x_0x_2 \\
 x_1x_2-x_0x_3\\
  x_2^2-x_1x_3\\
  
   x_1x_3^{v-1}-x_0^{1+\mu}+Ux_0x_4^{\mu n_0} \\
   x_2 x_3^{v-1}-x_0^{\mu}x_1+Ux_1x_4^{\mu n_0} \\
   x_{3}^v-x_0^{\mu}x_{2}+Ux_2x_4^{\mu n_0} \end{array}\right)$\\
   
  define the required deformation $\pi:Y\longrightarrow \mathbb{A}^1_\mathbb{ F}$. \item \ If \  $b=2$, \ $F=\left(\begin{array}{lllll}
  x_1^2-x_0x_2+Ux_0 x_4^{n_2}\\
 x_1x_2-x_0x_3+Ux_1x_4^{n_2}\\
  x_2^2-x_1x_3-U^2x_4^{2 n_2 }\\
   x_2x_3^{v-1}-x_0^{1+\mu}+Ux_2  x_4^{(v-1) n_3}+Ux_3^{v-1}x_4^{n_2}+U^2x_4^{(v-1) n_3+n_2}\\

     x_{3}^v-x_0^{\mu}x_{1}+Ux_3x_4^{(v-1) n_3}\end{array}\right)=$
     $$=f+U\left(\begin{array}{lllll}
   x_0 x_4^{n_2}\\
  x_1x_4^{n_2}\\
0\\
    x_2  x_4^{(v-1) n_3}+ x_3^{v-1}x_4^{n_2}\\
    x_3x_4^{(v-1) n_3}\end{array}\right)+U^2\left(\begin{array}{lllll}
  0\\
 0\\
  x_4^{2 n_2 }\\
  x_4^{(v-1) n_3+n_2}\\
    0\end{array}\right) \quad$$  
  define the required deformation $\pi:Y\longrightarrow \mathbb{A}^1_\mathbb{ F}$.   \denu
 \end{thm}
{\it Proof.} ({\bf Case $b=0$}).     In this case the image  of the element found in  \ref{as4} for  $T^1_B(-\mu n_0)$ (eigenvector) and the relation matrix among the generators of the ideal are: \ $g_1=(0,0,0,1)^T$,
$$   r=\left(\begin{array}{ccccccrrr}
x_2 &-x_1&x_0 &0  \\
-x_3 & x_2&-x_1 &0\\
x_0^{\mu}-x_3^{v} &0&0 &x_1^2-x_0 x_2 \\
0&x_0^{\mu}-x_3^{v}&0&x_1x_2-x_0x_3\\
0&0&x_0^{\mu}-x_3^{v} &x_2^2-x_1x_3\\

\end{array}\right)=\left(\begin{array}{ccccccrrr}
x_2 &-x_1&x_0 &0  \\
-x_3 & x_2&-x_1 &0\\
-f_4 &0&0 &f_1 \\
0&-f_4&0&f_2\\
0&0&-f_4 &f_3\\

\end{array}\right).$$
Hence an infinitesimal deformation of $X$ is given by the equations $F=f+\varepsilon Ug_1$.\\ 
 By a direct computation one has 
$$rUg_1=U \left(\begin{array}{cc} 0\\ 0\\  f_1\\ f_2\\f_3 \end{array}\right), \quad \rho=\left(\begin{array}{rrrcccrrr}
0&0 & 0 &0 \\
0&0 & 0 &0  \\
-U&0 & 0 &0  \\
0&-U & 0 &0  \\
0&0 & -U &0  \\
 \end{array}\right),\quad \rho g_1=\left(\begin{array}{cc} 0\\ 0\\ 0\\ 0\\0\\ 0\\ 0 \end{array}\right). $$ 
 Following the Pinkham's method \cite[(1.16)]{pi}, we consider the weighted homogeneous projective space $Proj(\mathbb{F}[x_0,\dots,x_4]),\quad \textrm{weight}(x_i)=n_i,$ for \ $0\leq i\leq 3, \  \textrm{weight}(x_4)=1)$ and substitute the variable  $U $   with $Ux_4^{ \mu  n_0},\  (U$ parameter); therefore  we get the   deformation $Y$ with parameter space $S=\mathbb{A}^1$ and fibres which are projective curves with only one (regular) point at infinity $P_{\infty}=(t^{n_0}:\dots:t^{n_3}:0), \ t\neq 0$. \  The equations are
  $$F=f+Ug_1x_4^{\mu n_0}.$$ To verify that the fibres $Y_U, \ U\neq 0$ of the family are  non-singular curves  it suffices to put $x_4=1$ and study the rank of the 
    jacobian matrix of the affine curve $Y_U\cap(x_4\neq 0$). Now this  matrix is equal to the jacobian matrix $J$ of the curve $X$.  We already know that   $rank_P(J) = 3$ if $P\neq (0,\dots,0)$. $$J =\left(\begin{array}{ccccccccc}
- x_2&2x_1 &-x_0\ \ &0\\
- x_3& x_2&x_1 &-x_0 \\
0&- x_3&2x_2 &-x_1 \\
-\mu x_0^{\mu-1}&0&0&vx_3^{v-1}\\
\end{array}\right).$$
Since by the equations of $Y_U$ we get $P\in Y_U\I P\neq (0,\dots,0)$, we are done.\\

\noindent({\bf Case $b=1$}).     In this case the image  of the element found in  \ref{as4} for  $T^1_B(-\mu n_0)$ (eigenvector) and the relation matrix among the generators of the ideal are 
$$   g_1:=\left(\begin{array}{cc} 0\\ 0\\ 0\\ x_0\\ x_1\\x_2 \end{array}\right),
 \quad r=\left(\begin{array}{rrrcccrrr}
x_2 &-x_1&x_0 &0 & 0 &0\\
-x_3 & x_2&-x_1 &0&0&0\\
x_3^{v-1}&0&0 &-x_1&x_0&0 \\
0&x_3^{v-1}&0&-x_2&0& x_ 0 \\
x_ 0^{\mu}&-x_3^{v-1}&0&0&x_1&-x_ 0\\
0& x_ 0^{\mu}&-x_3^{v-1}&-x_3&x_2&0\\
0&0&-x_3^{v-1}&0&x_2&-x_1\\
0&0&x_0^{\mu}&0&-x_3&x_2
\end{array}\right).$$
By a direct computation one has 
$$rUg_1=U \left(\begin{array}{cc} 0\\ 0\\ 0\\ 0\\f_1\\ f_2\\ 0\\f_3 \end{array}\right), \quad \rho=\left(\begin{array}{rrrcccrrr}
0&0 & 0 &0&0 & 0 \\

0&0 & 0 &0&0 & 0 \\
0&0 & 0 &0&0 & 0 \\
0&0 & 0 &0&0 & 0 \\
-U&0 & 0 &0&0 & 0 \\
0&-U & 0 &0&0 & 0 \\
0&0 & 0 &0&0 & 0 \\
0&0 & -U &0&0 & 0 \\
 \end{array}\right),\quad \rho g_1=\left(\begin{array}{cc} 0\\ 0\\ 0\\ 0\\0\\ 0\\ 0 \end{array}\right). $$ 
 Following the Pinkham's method as above, we substitute the variable  $U $   with $Ux_4^{ \mu  n_0},\  (U$ parameter); therefore  we get the   deformation $Y$ with parameter space $S=\mathbb{A}^1$ and fibres which are projective curves with only one regular point at infinity $P_{\infty}=(t^{n_0}:\dots:t^{n_3}:0), \ t\neq 0$.\\ The equations are
  $$F=f+Ug_1x_4^{\mu n_0}.$$ To verify that the fibres $Y_U, \ U\neq 0$ are  non-singular curves  it suffices to put $x_4=1$ and study the rank of the 
    jacobian matrix of the affine curve $Y_U\cap(x_4\neq 0)$. This matrix is
$$J =\left(\begin{array}{ccccccccc}
- x_2&2x_1 &-x_0\ \ &0\\
- x_3& x_2&x_1 &-x_0 \\
0&- x_3&2x_2 &-x_1 \\
-(1+\mu)x_0^{ \mu}+U& x_3^{v-1}&0&(v-1)x_1x_3^{v-2}\\
-\mu x_0^{\mu-1}x_1&-x_0^{\mu}+U& x_3^{v-1}&(v-1) x_2x_3^{v-2}\\
-\mu x_0^{\mu-1}x_2&0&-x_0^{\mu}+U&vx_3^{v-1 } \\
\end{array}\right).$$
We claim that the rank of $J =3 $ if $U\neq 0$ hence by the jacobian criterion of regularity we deduce that the fibres ar smooth for $U\neq 0$, i.e.,  the curve $X$ is smoothable and the semigroup is Weierstrass.\\
In $P_0(0:\dots:0,1),\ U\neq 0$ we have the non-null minor 
$\quad \textrm{det}\left(\begin{array}{ccccccccc}

 U& 0&0\\
0&U& 0\\
0&0&U\\
\end{array}\right).$\\
If $P\neq P_0$ ($P$ belonging to the fibre $Y_U$ of the canonical projection $\pi:Y\longrightarrow \mathbb{A}^1$),  according to  the equations of $Y_U$ \ and since \  $v\geq 2,$ \ we have 
$$\textrm{det}\left(\begin{array}{ccccccccc}
 2x_1 &-x_0\ \ &0\\
  x_2&x_1 &-x_0 \\
 
  \ \ x_3^{v-1}&0&(v-1)x_1x_3^{v-2}\\
\end{array}\right)= x_3^{v-1}x_0^2+(v-1)x_1x_3^{v-2}(2x_1^2+x_0x_2)=$$$=x_3^{v-2}x_0[x_0x_3+3(v-1)x_1x_2]=(3v-2)x_0^2x_3^{v-2}.$\\
If $x_0 x_3=0$,    by the equations we get only $P_0$. We are done.\\
 
  ({\bf Case $b=2$}). In this case one can easily see that   the generator found in \ref{as4} gives a deformation with all singular fibres.
 Then we need to find a different suitable deformation $Y$  such that  the rank of the jacobian matrix is "generically"  equal to 3  by the Jacobian criterion of regularity. We claim that a deformation which verifies this condition  is 
 $$F=f+Ug_1+Vg_2+ h\hspace{1.5cm}(*)$$
 where $g_1$,  $g_2$ are the images of the basis of $\ T^1_B(-(v-1)n_3)$ \ (resp $T^1_B(-n_2)$)  of  (\ref{as4}.3) and $h$ is found by the flatness conditions.\ \quad \
   Precisely, \ the relation matrix $\ r\ $ among the generators of the ideal and the vectors $g_1, g_2, h$ \  are  
$$r=\left(\begin{array}{cccccccc}
x_2 &-x_1\ &x_0 &0 & 0 \ \\
-x_3 & \ x_2&-x_1 &0&0 \ \\
0&  x_3^{v-1}&0& -x_1&x_0 \  \\
 x_ 0^{\mu}&0  &x_3^{v-1}& \  -x_2&x_ 1\\
0& x_ 0^{\mu}&0& -x_3&x_2 \\
\end{array}\right), \ \   g_1 =\left(\begin{array}{cc} 0\\ 0\\ 0\\ x_2\\ x_3 \end{array}\right), \ \ g_2=\left(\begin{array}{cc} x_0\\ x_1\\ 0\\ \ x_3^{v-1}\\0\\ \end{array}\right),\ \ h=\left(\begin{array}{cc}    0\\ 0\\ -V^2\ \ \\ UV\\  0 \end{array}\right).$$
To find $h$, by a direct computation we get   $$rUg_1=\left(\begin{array}{cc} 0\\ 0\\ -f_2\\ -f_3\\ 0\end{array}\right)=-\rho_1 f,\ \ {\rm with\ \ } \rho_1= \left(\begin{array}{cccccccc}
0&0&0 &0 & 0  \\
0&0&0 &0 & 0  \\
0&U&0 &0 & 0  \\
0&0&U &0 & 0  \\
0&0&0 &0 & 0  \\

\end{array}\right),$$
  $$rVg_2=V\left(\begin{array}{r} -f_1\\ f_2\\ 0\ \ \\ -f_4\\ -f_5\end{array}\right)=-\rho_2 f,\ \ {\rm with\ \ } \rho_2= \left(\begin{array}{cccccccc}
V&0&0 &0 & 0  \\
0&-V\ &0 &0 & 0  \\
0&0&0 &0 & 0  \\
0&0&0 &V & 0  \\
0&0&0 &0 & V  \\
\end{array}\right), $$
Then $\rho = \left(\begin{array}{cccccccc}
V&0&0 &0 & 0  \\
0&-V \ \ &0 &0 & 0  \\
0&U&0 &0 & 0  \\
0&0&U &V & 0  \\
0&0&0 &0 & V  \\
\end{array}\right),\ \ \rho g=\left(\begin{array}{cc}\ \ V^2x_0\\ -V^2x_1\\ x_1 UV\\ x_2UV+V^2x_3^{v-1}\\  x_3UV \end{array}\right)$.\\
Further \ \ $\rho g=-r h$, \ \ with \ \ $h=\left(\begin{array}{cc}    0\\ 0\\ -V^2\ \ \\ UV\\  0 \end{array}\right)$, \ \ finally \ \ $\rho  h=\left(\begin{array}{cc} 0\\ 0\\ 0\\ 0\\  0 \end{array}\right).$\\
Following the Pinkham's method as above, we substitute the variables $U,V$ respectively with $Ux_4^{(v-1) n_3}, Vx_4^{n_2}$; therefore  we get the   deformation $Y'$ with parameter space $S=\mathbb{A}^2$ and fibres    with only one (regular) point at infinity. The equations are $$\left(\begin{array}{lllll}
  x_1^2-x_0x_2+Vx_0 x_4^{n_2}\\
 x_1x_2-x_0x_3+Vx_1x_4^{n_2}\\
  x_2^2-x_1x_3-V^2x_4^{2(n_2)}\\
   x_2x_3^{v-1}-x_0^{1+\mu}+Ux_2  x_4^{(v-1) n_3}+Vx_3^{v-1}x_4^{n_2}+UVx_4^{(v-1) n_3+n_2}\\
     x_{3}^v-x_0^{\mu}x_{1}+Ux_3x_4^{(v-1) n_3}\end{array}\right).$$
Finally we claim that the restriction to the line $(U=V)\subseteq S$ gives a 1-parameter deformation $Y$ with smooth generic fibre. Since the point $P_{\infty}$ is   non-singular,   we can put $x_4=1$ and   study  the rank of 
  jacobian matrix of the affine curve $Y\cap (x_4\neq 0)$:   
$$J(Y)=\left(\begin{array}{ccccccccc}
- x_2+U&2x_1 &-x_0 &0\\
- x_3& x_2+U&x_1 &-x_0 \\
0&- x_3&2x_2 &-x_1 \\
-(1+\mu)x_0^{ \mu} & 0&x_3^{v-1}+U&(v-1)x_3^{v-2}(x_2+U)\\
-\mu x_0^{\mu-1}x_1&-x_0^{\mu} &0& v x_3^{v-1}+U\\
 \end{array}\right).$$
 As in case $b=1$ we can assume $P\neq P_0=(0:\dots : 0:1)$, $P$ belonging to the fibre $Y_U$. \ Consider the minor
$$\textrm{det}\left(\begin{array}{ccccccccc}
 2x_1 &-x_0 &0\\
  x_2+U&x_1 &-x_0 \\
 - x_3&2x_2 &-x_1 \\
  \end{array}\right)=\textrm{det}\left(\begin{array}{ccccccccc}
 2x_1 &-x_0 &0\\
  x_2 &x_1 &-x_0 \\
 - x_3&2x_2 &-x_1 \\
  \end{array}\right)+\textrm{det}\left(\begin{array}{ccccccccc}
 2x_1 &-x_0 &0\\
   U&x_1 &-x_0 \\
 - x_3&2x_2 &-x_1 \\
  \end{array}\right)=$$
  $=-Ux_0 x_1$. \vspace{0.2cm}\\
   If $x_0=0$ by the equations we get      $x_1 =0$,  $ x_3^{v-1}+U =0  $ and $x_2=+ U$, and so by the fourth equation give $Ux_3^{v-1}=0$, impossible. \\
     If $x_1=0$ by the equations we get $x_2=U$ (since $x_0\neq 0$ by above),\ $x_3=0$.\\  The fourth equation  gives \ \ $     -x_0^{\mu+1}+U^2=0$. \\ Now the jacobian matrix in these points is
$$
J(Y)=\left(\begin{array}{ccccccccc}
0&0 &-x_0 &0\\
0& 2U&0 &-x_0 \\
0&0 &2U &0 \\
-(1+\mu)x_0^{ \mu} & 0& U&(v-1)x_3^{v-2}(x_2+U)\\
0&-x_0^{\mu} &0&  U\\
 \end{array}\right).
 $$
 Recalling the fourth equation we see that the minor 
 $$\textrm{det}\left(\begin{array}{ccccccccc}
  2U&0 &-x_0 \\
 0 &2U &0 \\
 -x_0^{\mu} &0&  U\\
 \end{array}\right)= 2U(2U^2-x_0^{\mu+1})=2U^3
 $$ 
 Hence $\textrm{rank}(J(Y))=3$ for each $P\in Y_U,\ \forall \ U\neq 0$.
  We are done.

  \begin{rem}
  {\rm Note that in case   $b=2, v>2$  the deformation $Y'$ of the curve $X$ with equations   $$\left(\begin{array}{lllll}
  x_1^2-x_0x_2 \\
 x_1x_2-x_0x_3\\
  x_2^2-x_1x_3\\
   x_1x_3^{v-1}-x_0^{1+\mu}+Vx_0\\
   
   x_{3}^v-x_0^{\mu}x_{1}+Vx_1\end{array}\right)$$  has parameter space $\mathbb{A}^1_\mathbb{ F}$, but
       every fibre of this deformation has a singularity at the origin: hence in general this  construction  does not give informations on the smoothability of the curve $X$ even if the algorithm  to construct $Y'$ starting from the infinitesimal deformation ends at the first step. }  
  \end{rem}

\makeatletter
\newcommand\ackname{Acknowledgement}
\if@titlepage
  \newenvironment{acknowledgements}{%
      \titlepage
      \null\vfil
      \@beginparpenalty\@lowpenalty
      \begin{center}%
        \bfseries \ackname
        \@endparpenalty\@M
      \end{center}}%
     {\par\vfil\null\endtitlepage}
\else
  \newenvironment{acknowledgements}{%
      \if@twocolumn
        \section*{\abstractname}%
      \else
        \small
        \begin{center}%
          {\bfseries \ackname\vspace{-.5em}\vspace{\z@}}%
        \end{center}%
        \quotation
      \fi}
      {\if@twocolumn\else\endquotation\fi}
\fi
\makeatother

{\footnotesize
\baselineskip11pt

\renewcommand{\refname}{{\large\bf References}}

\end{document}